\documentclass{amsart}
\usepackage[colorlinks=true,pagebackref,hyperindex,citecolor=blue,linkcolor=blue]{hyperref}
\usepackage[top=1in, bottom=1in, left=1.1in, right=1.1in]{geometry}
\usepackage{amsmath}
\usepackage{amsfonts}
\usepackage{amssymb}
\usepackage{amsthm}
\usepackage{stmaryrd}
\usepackage{braket}
\usepackage[all]{xy}
\usepackage{mathrsfs}
\usepackage{enumitem} 
\usepackage[T1]{fontenc}
\usepackage{color}
\usepackage{tikz}
\usepackage{nicefrac, xfrac}
\usepackage[numbers]{natbib}
\usepackage{comment} 

\newcommand*\rfrac[2]{{}^{#1}\!/_{#2}}

\newtheorem{theorem}{Theorem}[section]
\newtheorem{corollary}[theorem]{Corollary}

\newtheorem{lemma}[theorem]{Lemma}
\newtheorem{proposition}[theorem]{Proposition}

\newtheorem*{theoremx}{Theorem}

\theoremstyle{definition}

\numberwithin{equation}{section}

\newcommand{\ZZ}{\mathbb{Z}}
\newcommand{\PP}{\mathbb{P}}

\newcommand{\CC}{\mathbb{C}}

\newcommand{\ctM}{\Theta_{M}}
\newcommand{\mcL}{\mathcal{L}}
\newcommand{\mcO}{\mathcal{O}}

\newcommand{\mcH}{\mathcal{H}}

\newcommand{\mcE}{\mathcal{E}}
\newcommand{\mcF}{\mathscr{F}}

\makeatletter
\def\@tocline#1#2#3#4#5#6#7{\relax
  \ifnum #1>\c@tocdepth
  \else
    \par \addpenalty\@secpenalty\addvspace{#2}%
    \begingroup \hyphenpenalty\@M
    \@ifempty{#4}{%
      \@tempdima\csname r@tocindent\number#1\endcsname\relax
    }{%
      \@tempdima#4\relax
    }%
    \parindent\z@ \leftskip#3\relax \advance\leftskip\@tempdima\relax
    \rightskip\@pnumwidth plus4em \parfillskip-\@pnumwidth
    #5\leavevmode\hskip-\@tempdima
      \ifcase #1
       \or\or \hskip 1.9em \or \hskip 2em \else \hskip 3em \fi%
      #6\nobreak\relax
    \dotfill\hbox to\@pnumwidth{\@tocpagenum{#7}}\par
    \nobreak
    \endgroup
  \fi}
\makeatother

\begin{document}

\title[Foliations on K3 Surfaces]{Foliations on Projective Complete Intersection K3 Surfaces}

\author[J. Olivares]{Jorge Olivares}
\address{Centro de Investigaci\'on en Matem\'aticas, A.C.
A.P. 402, Guanajuato 36000, México} 
\email{olivares@cimat.mx}

\author[D. Posada]{Daniel Posada-Buriticá{$^1$}}
\address{Centro de Investigaci\'on en Matem\'aticas, A.C.
A.P. 402, Guanajuato 36000, México} 
\email{daniel.posada@cimat.mx}

\thanks{{$^1$}This author was partially supported by CONACYT-SECIHTI Scholarship 842809.}



\begin{abstract}
We study foliations $\mathscr{F}$ on projective complete intersection K3 surfaces $X \hookrightarrow \PP^n$, where $\mathscr{F}$ has isolated singularities and it is the restriction of a foliation of degree $ d $ on $ \PP^n $ that leaves $ X $ invariant. We compute the values of the degrees $ d $ for which $ \mathscr{F} $ is uniquely determined by its singular scheme.
\end{abstract}

\maketitle

\section*{Introduction}
Let $ M $ be a compact connected complex manifold of dimension $ n $. Recall that a
holomorphic foliation by curves $ \mathscr{F} $ on $ M $  (simply \emph{a foliation on} $ M $
in the sequel) may be defined by
non-identically zero holomorphic vector fields $ X_i $ defined on a covering $ \{ V_i \} $ of $ M $
such that in each overlapping set $ V_i \cap V_j \neq \emptyset $ we have
\begin{equation}\label{cocicle}
   X_i = \xi_{ij} X_j,
\end{equation}
where $ \xi_{ij} $ is a never vanishing holomorphic function. If
$ L^{\vee} $ denotes the holomorphic line bundle constructed with the cocycle $( \xi_{ij} )$,
and $\mcL^{\vee}$ its corresponding invertible sheaf,
then the $ X_i $'s give rise to a global section
$ s \in H^{0}(  M , \ctM \otimes \mcL^{\vee}) $ or to a global section in
$ H^{0}( M, \mathscr{H}om_{ \mcO_M} (\mcL , \ctM) )$, where
$ \ctM $ is the tangent sheaf of $ M $ and $\mcL$ is the dual
of $\mcL^{\vee}$. Two global sections (in the corresponding spaces) define the same
foliation if and only if one is a non-zero scalar multiple of the other.

The line bundle
$ L^{\vee} $ just defined 
will be called the \emph{cotangent bundle of} $\mathscr{F}$ and its dual $ L $, \emph{its tangent bundle}. Hence, the space $\mathscr{F}ol(\mcL, M) $ of foliations
$ \mcF $ with tangent bundle $ L $ (or tangent sheaf $ \mcL $) is $ \PP H^{0}( M, \mathscr{H}om(\mcL , \ctM) )$. Such an $ \mcF $ corresponds to a foliation with cotangent bundle $ L^{\vee}$ (or cotangent sheaf 
$ \mcL^{\vee} $)
by regarding it as the class $ [ s ] \in \PP H^{0}( M, \ctM \otimes \mcL^{\vee})$ of a global section $ s \in H^{0}( M, \ctM \otimes \mcL^{\vee}) $.

Given a global section
$ s \in H^{0}( M, \mathscr{H}om_{ \mcO_M} (\mcL , \ctM) )$,
the scheme $ Z = Z_s $ of those points $ p\in M $ where the induced morphism
$\mcL_p\rightarrow \Theta_{M,p}$ is zero will be referred to as
the \emph{singular scheme} of $ s $: its ideal sheaf 
$ \mathcal{I}_Z \subset \mcO_M $
is the sheaf obtained by gluing the ideals $ (a_i^1, \ldots, a_i^n) \subset \mcO( V_i ) $,
where $ a_i^1,\ldots, a_i^n $ are the components of the vector field $ X_i $ that defines
$ s $ on the open set $ V_i $, as described nearby \eqref{cocicle}. It is clear that the singular scheme
of $ \mcF = [s] $ is the singular scheme of any section in $ [s] $.

We say that $ \mcF = [s]$ has \emph{isolated singularities} if
$ \text{dim } Z = 0 $.

For foliations on projective spaces $ \PP^n $, with $ n \geq 2 $, we have the following result
from \cite{CampilloOlivares}: for $ d \geq 2 $, a foliation
$ [s] \in \PP H^0( \PP^n, \mathscr{H}om( \mcO_{\PP^n}( 1-d ), \Theta_{\PP^n} ) )
 = \mathscr{F}ol( \mcO_{\PP^n}( 1-d ), \PP^n ) $
 with isolated singularities 
\textit{is the unique one in} $\mathscr{F}ol ( \mcO_{\PP^n}( 1-d ), \PP^n ) $
 \textit{with singular scheme} $ Z = Z_s $.

 We summarize this statement saying that, for $ d \geq 2 $,
 a foliation with isolated singularities 
 and tangent sheaf 
  $ \mcO_{\PP^n}( 1-d ) $
 on a projective space of dimension $ n \geq 2 $ is uniquely determined by its singular scheme.

 In this paper we study the extension of this result to foliations with isolated singularities on complete intersection projective K3 surfaces $ X $.

Other results of this type, dealing with foliations (or distributions) of rank 
equal to $1$ or higher
on projective spaces, are given in \cite{G-P}, \cite{C-F-N-V} or in \cite{A-C}: 
the first and second consider, respectively, codimension one 
 and dimension one foliations with $1$-dimensional singular scheme, both in $ \PP^3 $;
 the third one, which is a far reaching generalization of the main result in 
 \cite{CampilloOlivares}, considers higher rank distributions with $1$-codimensional or 
 with $0$-dimensional singular scheme. Each of them exhibits sufficient conditions on the
 singular scheme $ Z $ of the foliation or distribution $ \mathscr{F} $
 in order that $ \mathscr{F} $ becomes uniquely determined by $ Z $. The case of
 foliations on 
 Hirzebruch surfaces $ S_{\delta}, \, \delta \geq 0 $, is treated in \cite{Hirzebruch}.

 The quoted result from \cite{CampilloOlivares} traces back to \cite{GMKempf}, where it was
 proved for foliations with \textit{reduced} singular scheme. 
 This hypothesis was first removed for foliations on the projective plane in 
 \cite{Polarity}, and this result has had applications in the problem of classification of 
 foliations (see for instance \cite{Claudia}, \cite{ClaudiaRubi} or \cite{ClaudiaRonzon}).
 The results presented here could have similar applications.

Now, to state our results, let
$ X \overset{ i }{\hookrightarrow} \PP^n $
be a complete intersection projective K3 surface.
All along the paper, we will consider foliations on $ X $ with tangent sheaf 
$ \mathcal{ L }_d = i^{\ast}\mathcal{O}_{\PP^n}(1-d) $, 
where the dimension $ n $ depends on the complete intersection type of $ X $.
This means that we shall consider foliations 
$ \mathscr{F} = [ s ]$ on $ \PP^n $, with 
$ s \in H^0( \PP^n, \mathscr{H}om( \mcO_{\PP^n}( 1-d ), \Theta_{\PP^n} ) ) $,
such that its restriction $ i^{\ast} s $ to $ X $ factors through the tangent sheaf
$ \Theta_X $ of $ X $: 
$ i^{\ast} s \in H^0( X, \mathscr{H}om( i^{\ast}\mathcal{O}_{\PP^n}(1-d), \Theta_X ) ) $.
That is, foliations $ \mathscr{F} $ which leave $ X $ invariant.
As quoted in \cite{MarcioAnn}, the restriction $ i^{\ast} s $ may have
isolated singularities for an $ s $ with non-isolated ones.

Our results are summarized in the following statement:
\begin{theoremx}
Let $ X \overset{ i }{\hookrightarrow} \PP^n $ be a complete intersection projective K3 surface, 
and consider the invertible sheaves $\mathcal{L}_d=i^{\ast}\mathcal{O}_{\PP^n}(1-d)$ on $X$, with 
$ d \geq 0 $. Then,
\begin{itemize}
\item[ \textbf{A} ] \label{Theorem A}
Foliations on $ X $ with tangent sheaf $ \mcL_d $ exist only for $ d \geq 3 $.
\item[ \textbf{B} ]  \label{Theorem B}
Let $ d \geq 3 $. Then any foliation on $ X $ with isolated singularities and tangent sheaf 
$ \mathcal{L}_d $ is uniquely determined by its singular scheme for:
    \begin{itemize}
     \item[ (1) ] $d \geq 6$, in case $ X $ is a quartic hypersurface in $\PP^3$;
     \item[ (2) ] $d \geq 5$, in case $ X $ is a complete intersection of
     a quadric and a cubic in $\PP^4$, and 
     \item[ (3) ] $d \geq 4$, in case $ X $ is a complete intersection of
     three quadrics in $\PP^5$.
    \end{itemize}     
\end{itemize}
\end{theoremx}
 Every complete intersection projective K3 surface $ X $ fits into one of the cases listed in Theorem B above (see Lemma \ref{K3-CI} below). 
 In each case, the Picard rank of $ X $ may be strictly greater
 than $ 1 $ but it is always $ 5\leq 20 $ (\cite[Chapter 1, $\mathsection 3.3$]{Huybrechts}). However,
 the \textit{very general} $ X $ has Picard rank
 equal to $ 1 $ (\cite[Chapter 16, $\mathsection 4.4$]{Huybrechts}) and hence, in these cases, our 
 theorem describes all possible foliations on $ X $. 
 
 Every foliation as in Theorem A above has a nonempty singular scheme:
 see Proposition \ref{IsoSings} below (which hence has no content
 for $ d= 0, 1 $ and $ 2 $).

 Recall that Theorem A is the solution of the so-called Poincaré's Problem 
 for these surfaces $ X $: Case (1) is a straightforward consequence of \cite{MarcioInv}
 and, in Case (2), we have obtained the same conclusion as that from 
 \cite[Theorem II]{MarcioAnn},
 although this result applies only for complete intersections of odd dimension.

 \vspace{5mm}

The structure of the paper is the following: although somehow standard, the notation just used will be presented
in Section \ref{Bckgnd}, where we will also review the results in K3 surfaces,
complete intersections and foliations that will be used along the paper.
Both proofs of Theorem A and B are divided into the three types of K3 surfaces we are dealing with:
The statement in Theorem A above is a short version of the computation of the dimensions
$ h^0(X, \Theta_X \otimes i^{\ast}\mathcal{O}_{\PP^n}(d-1)) $
of the spaces of sections associated to foliations with tangent sheaf 
$ \mathcal{ L }_d = i^{\ast}\mathcal{O}_{\PP^n}(1-d) $:
these will be presented in
Section \ref{Ch:PffThmA} (Theorems 
\ref{Thm:Ex-1}, \ref{Thm:Ex-2} and 
\ref{Thm:Ex-3} below). The main idea for the proof of Theorem B is 
presented
at the beginning of Section \ref{Ch:PffThmB}, where 
the proof itself will be given 
(see Theorems \ref{Thm:Uniq-1}, 
\ref{Thm:Uniq-2} and \ref{Thm:Uniq-3} below).

\section{Preliminaries}\label{Bckgnd}
 Let $ \PP^n = \textup{Proj} \ \CC[ x_0, \dots, x_n ] $ be the complex projective space of dimension 
 $ n \geq 2 $, and let $ \mcO_{\PP^n}$, $ \Theta_{\PP^n} $, $ \Omega^1_{\PP^n} $ and $ \mathcal{H} $
 denote its structure, tangent, cotangent and hyperplane sheaves, respectively. For an 
 $ \mcO_{\PP^n}$-sheaf $ \mathcal{E} $, we will write 
 $ \mcE( d ) $ for $ \mcE \otimes \mcH^{\otimes d}$, if $ d \geq 0 $, and 
$ \mcE \otimes (\mcH^{\vee})^{\otimes |d|}$, if $ d < 0 $. Similar notation will be adopted
for a submanifold (or subvariety) $ X \subset \PP^n $. Any vector bundle
will be identified with its locally free sheaf of sections.

\subsection{Projective K3 Surfaces and complete intersections:} Here we recall some facts on K3 
surfaces (mostly following the comprehensive monograph \cite[]{Huybrechts}), and on projective
complete intersections, that will be used in the sequel.

An \emph{algebraic K3 surface} over $\CC$  is a complete nonsingular variety $X$ over $\CC$ of dimension $2$, whose canonical sheaf $ \omega_X = \wedge^2 \Omega_X^1 $ is isomorphic to 
$\mathcal{O}_X$ and with $H^1(X, \mathcal{O}_X)=0$. Since every nonsingular complete surface over 
$\CC$ is projective (see \cite[Remark 4.10.2 (b)]{Hartshorne}), algebraic K3 surfaces are always projective.

A \emph{complex K3 surface} is a compact connected complex manifold $X$ of dimension $2$ with 
$ \omega_X \simeq \mathcal{O}_X $
and $H^1(X, \mathcal{O}_X)=0$. In particular, the triviality of $\omega_X$ implies the existence of a non-zero holomorphic $2$-form $\varphi \in H^0(X, \omega_X)$. \emph{Contracting a vector 
field with $\varphi$ induces a non-canonical isomorphism $\Theta_X \cong \Omega^1_X$}.

Since $\PP^n$ carries both an algebraic and a complex structure, the GAGA principle 
\cite[Appendix B]{Hartshorne} yields a natural embedding
 \[ \{\textup{Algebraic K3 surfaces over $\CC$}\} \hookrightarrow \{\textup{Complex K3 surfaces}\},\]
whose image is the set of complex K3 surfaces that are projective  
(see \cite[Chapter I, $\mathsection 3.1$]{Huybrechts}). This correspondence allow us to use freely 
both the algebraic and complex structures on such surfaces. \\

Let $X=V(f_1, \ldots, f_c) \subset \PP^n$ be a (smooth) complete intersection of codimension $ c $, where each $f_i$ is a homogeneous polynomial of degree $d_i$. We shall refer to such an $X$ as a 
\emph{(smooth) complete intersection of type $(d_1, \ldots, d_c)$ in $ \PP^n $}. We assume that 
$ d_i > 1 $ for every $i$ to guarantee that $ n $ is minimal.

The augmented Koszul complex \cite[Chapter IV, $\mathsection 2$]{RRAlgebra} 
associated to such an $ X $
gives the following resolution of the \emph{ideal sheaf} $\mathcal{I}_X$ of $X$:
\begin{equation}\label{KoszCxSh}
0 \longrightarrow \mathcal{O}_{\PP^n}\left(-\sum_{i=1}^c d_i \right) \longrightarrow  \cdots \longrightarrow \bigoplus_{i \neq j} \mathcal{O}_{\PP^n}(-d_i-d_j) \longrightarrow \bigoplus_{i=1}^c \mathcal{O}_{\PP^n}(-d_i) \longrightarrow \mathcal{I}_X \longrightarrow 0.
\end{equation}

It is not hard to see from \eqref{KoszCxSh} that
the \emph{normal sheaf} $\mathcal{N}_{X/\PP^n}$ of $X$ in $\PP^n$, defined as the dual of the conormal sheaf $\mathcal{I}_X/ \mathcal{I}_X^2$, is given by
\begin{equation}\label{NormalBundle}
\mathcal{N}_{X/\PP^n} \cong \bigoplus_{i=1}^ci^{\ast} \mathcal{O}_{\PP^n}(d_i).
\end{equation}
The normal sheaf fits into the exact sequence \eqref{NormalSeq} below. 

For the following result, see \cite[\S 2.2]{Cr-Es}:
\begin{lemma}\label{aCM}
Let $ X \subset \PP^n $ be a complete intersection with ideal sheaf $\mathcal{I}_X$. Then
$ H^j( \PP^n, \mathcal{I}_X( m ) ) = 0 $ for every $ m \in \ZZ $ and 
$ j = 1, \dots, \textup{dim}(X) $.
\end{lemma}
\begin{lemma}[{\cite[Example 1.3.]{Huybrechts}}] \label{K3-CI}
The only complete intersection projective K3 surfaces $ X $ are the following:
\begin{itemize}
    \item[(1) ] Smooth quartic hypersurfaces in $\PP^3$,
    \item[(2) ] Smooth complete intersections of type $(2,3)$ in $\PP^4$, and
    \item[(3) ] Smooth complete intersections of type $(2,2,2)$ in $\PP^5$.
\end{itemize}    
\end{lemma}
The proof follows from \cite[Exercise. 8.4]{Hartshorne}:
the canonical sheaf of $ X $ is given by $ \omega_X = \mcO_{X}( \sum d_i - n - 1 ) $: the
only possible solutions in $ n $ and in the $ d_i $ to the equation
$ \sum d_i = n + 1 $ are the ones in the list. In each case, one then
verifies that the condition $ H^1( X, \mcO_X ) = 0 $ holds.

For further reference, we state the following easy consequence of \eqref{KoszCxSh}:
\begin{corollary}\label{K3-CI-IdShf}
      Consider the three types of complete intersections $ X $ from Lemma \ref{K3-CI}. Then their ideal sheaves $\mathcal{I}_X$ admit respectively
      the following resolutions:
    \begin{itemize}
        \item[(1) ] $0 \longrightarrow \mathcal{O}_{\PP^3}(-4) \longrightarrow \mathcal{I}_X \longrightarrow 0$,
        \item[(2) ]  $0 \longrightarrow \mathcal{O}_{\PP^4}(-5) \longrightarrow 
        \mathcal{O}_{\PP^4}(-2)\oplus \mathcal{O}_{\PP^4}(-3) \longrightarrow 
        \mathcal{I}_X \longrightarrow 0$, and
        \item[(3) ] 
        $0 \longrightarrow \mathcal{O}_{\PP^5}(-6) \longrightarrow 
        \mathcal{O}_{\PP^5}(-4)^{\oplus 3}  \longrightarrow \mathcal{O}_{\PP^5}(-2)^{\oplus 3} \longrightarrow \mathcal{I}_X \longrightarrow 0$, 
        which can be written as the short exact sequence
      \;  $0 \longrightarrow 
    \left(\rfrac{\mathcal{O}_{\PP^5}(-4)^{\oplus 3}}{\mathcal{O}_{\PP^5}(-6)}\right) 
     \longrightarrow \mathcal{O}_{\PP^5}( -2 )^{\oplus 3} \longrightarrow  \mathcal{I}_X \longrightarrow 0.$
    \end{itemize}
\end{corollary}
\subsection{Foliations by Curves:}

Foliations by curves 
with tangent sheaf 
$ \mathcal{O}_{\PP^n}(1-d) $ on $\PP^n$ arise
from polynomial homogeneous vector fields of degree $ d $  in $ \CC^{n+1} $ through the following 
construction: Let $ \pi : \CC^{n+1} \setminus \{0\} \longrightarrow \PP^n ; \; p\mapsto [p] $, be the
projection map. Tensor the Euler sequence in $\PP^n$ \cite[p. 3]{Okonek} with the invertible sheaf
$\mathcal{O}_{\PP^n}(d)$ to obtain the exact sequence
\begin{equation*}
    0 \longrightarrow \mathcal{O}_{\PP^n}(d-1) \longrightarrow \mathcal{O}_{\PP^n}(d)^{\oplus (n+1)} \overset{ \pi_{*} }{\longrightarrow} \Theta_{\PP^n}(d-1) \longrightarrow 0,
\end{equation*}
whose long exact cohomology sequence is:
\begin{equation}\label{twisteulseq}
    0 \to H^0(\PP^n,\mathcal{O}_{\PP^n}(d-1)) \to H^0 (\PP^n,\mathcal{O}_{\PP^n}(d))^{\oplus (n+1)} 
    \overset{ \pi_{*} }{\longrightarrow} H^0(\PP^n,\Theta_{\PP^n}(d-1)) \to 0,
\end{equation}
since $H^1(\PP^n, \mathcal{O}_{\PP^n}(d-1))=0$ (see Lemma \ref{LemmBottFla} below).

It follows from \eqref{twisteulseq} that any section
$ s \in H^0( \PP^n, \Theta_{\PP^n}( d-1 ) ) \cong 
H^0( \PP^n, \mathscr{H}om( \mcO_{\PP^n}( 1-d ), \Theta_{\PP^n} ) ) $
 (and hence,
 any foliation on $\PP^n$ with tangent sheaf $ \mathcal{O}_{\PP^n}(1-d) $) 
 is defined by a polynomial homogeneous vector field of degree $d$ in $\CC^{n+1}$, where 
 two such vector fields $F$ and $G$ define the same section $ s $ iff $G= F + h \cdot R$, where
 $R=\sum_{i=0}^n x_i \frac{\partial}{\partial x_i}$ is the radial vector field in $\CC^{n+1}$ and
  $h$  is  some homogeneous polynomial of degree $d-1$. In particular, no such global sections
  $ s $ exist for $ d < 1 $.

 It follows moreover from \eqref{twisteulseq} that for $ 0 \neq s = \pi_{*} ( F ) $, with
 $F=\sum_{i=0}^n F_i( x ) \frac{\partial}{\partial x_i}$ as above, the singular \textit{set} of
 $ s $ consists in the points $ [p] \in \PP^n $ where $ R( p )\wedge F( p ) = 0 $
 and hence, the singular \textit{scheme}
 $ Z = Z_s $ of $ s $ is the subscheme of $\PP^n $ whose homogeneous ideal $ I_Z $
 is generated by the $ 2 \times 2 $ minors of the $ 2 \times (n+1) $ matrix with rows 
 $ x $ and $ F( x )$.

 These minors are the coefficients of the
 unique (up to a scalar multiple)
 projective $(n-1)$-form $ \omega_F $
 that annihilates 
$ s = \pi_{*} ( F ) $. One can compute it through the isomorphism
 \[ \Theta_{\PP^n}( d-1 ) \longrightarrow \wedge^{n-1} \Omega^1_{\PP^n}( n+d ) \quad ; \quad
     F \mapsto \omega_F :=\iota_{F}(\iota_{R}(dx_0 \wedge \cdots \wedge dx_n)),
   \]
 One can always assume that the codimension of $ Z $ is $ \geq 2 $,
 by requiring that the coefficients of $ F $ do not have a common factor.   

\subsubsection{Degree of the singular scheme of a foliation:}\label{DegofSS}
Let $ X $ be a compact complex surface and let
$ s \in H^0(X, \Theta_X \otimes \mathcal{L}^{\vee})$ be a section with isolated singularities,
locally defined by vector fields $ X_i $ as in \eqref{cocicle}, and let $ \mcF = [s] $.
Say 
$ X_i = a_i(x,y) \partial/\partial x + b_i(x,y) \partial/\partial y $ in $ V_i \subset X $.

For any point $ q \in V_i $, recall that the multiplicity 
$\mu_q( \mcF ) = \dim_{\CC}\left(\frac{ \mathcal{O}_{X, q} }{(a_i, b_i)\cdot 
 \mathcal{O}_{X, q} }\right)$
is a non-zero and finite positive integer only if $ q $ is an isolated
singular point of $ s $.

Now, following \cite[Proposition 2.4]{RuledSurfaces}, the number of singular points of $ s $
(counted with multiplicities) is equal to $ c_2(\Theta_X \otimes \mathcal{L}^{\vee})[ X ] $.
This number will be called the \textit{degree} of the singular scheme $ Z_s $ of $ s $. 
Now we proceed to compute them for the complete intersection K3 surfaces from 
Lemma \ref{K3-CI}:

By \cite[p. 28]{Friedman} we have
\begin{equation*}
  c_2(\Theta_X \otimes \mathcal{L}^{\vee})[X] =
    (c_2(\Theta_X)+c_1(\Theta_X) \cdot c_1(\mathcal{L}^{\vee})+c_1(\mathcal{L}^{\vee})^2)[X]\\
     = 24 + c_1(\mathcal{L}^{\vee})^2[X],
\end{equation*}
since $c_1(\Theta_X)[X] = c_1(X)[X] =0$ and $c_2(\Theta_X)[X]=c_2(X)[X] = 24 $, for any K3 surface 
$ X $, by \cite[p. 8]{Huybrechts}. Now we come to our case of interest: for
$ X \overset{ i }{\hookrightarrow} \PP^n $ and $\mathcal{L}=i^{\ast}\mathcal{O}_{\PP^n}(1-d)$,
we have
$ c_1(\mathcal{L}^{\vee})^2 = ( d - 1)^2 \cdot c_1(i^{\ast}\mathcal{O}_{\PP^n}(1))^2 $,
so that 
\[ c_1(\mathcal{L}^{\vee})^2[X] = ( d - 1 )^2 \cdot d_0, \]
where $ d_0 $ is the degree
of $ X \hookrightarrow \PP^n $, by \cite[p. 369]{MarcioAnn}. 
In sum, we have shown
\begin{proposition}\label{IsoSings}
 Let $ d \geq 0 $ be an integer, let
 $ X \overset{ i }{\hookrightarrow} \PP^n $ be a complete 
 intersection K3 surface
 and let $ s \in H^0( X, \Theta_X \otimes i^{\ast}\mathcal{O}_{\PP^n}(d-1) ) $ be a section
 with isolated singularities. Then the degree of the singular scheme $ Z_s $ of $ s $ is equal to
\begin{itemize}
    \item[1)] $ 24 + 4 \cdot ( d - 1 )^2 = 4\cdot( d^2-2d+7 ) $, 
        if $X$ is a smooth quartic hypersurface in $\PP^3$,
    \item[2)] $ 24 + 6 \cdot ( d - 1 )^2 = 6\cdot( d^2-2d+5 ) $, if $X$ 
     is of type $(2,3)$ in $\PP^4$,
    \item[3)] $ 24 + 8 \cdot ( d - 1 )^2 = 8\cdot( d^2-2d+4 )$, if $X$ 
     is of type $(2,2,2)$ in $\PP^5$.
    \end{itemize}
\end{proposition}
\section{Precise statements of Theorem A}\label{Ch:PffThmA}

Let $ d \geq 0 $ be an integer.
In this chapter we compute the dimensions
\[
h^0( X, \mathscr{H}om(i^{\ast}\mathcal{O}_{\PP^n}(1-d), \Theta_X )) =
 h^0( X, \Theta_X \otimes i^{\ast}\mathcal{O}_{\PP^n}(d-1) ) 
 = h^0( X,\Omega^1_X \otimes i^{\ast}\mathcal{O}_{\PP^n}(d-1) ),
\]
of the spaces of sections associated to foliations with tangent sheaf 
$ \mathcal{ L }_d = i^{\ast}\mathcal{O}_{\PP^n}(1-d) $,
for all projective complete intersection  K3 surfaces
$ X \overset{ i }{\hookrightarrow} \PP^n $
(Lemma \ref{K3-CI}).

The proofs of the theorems of this section follow a common pattern: Since $X$ is a nonsingular subvariety of $\PP^n$, we have the normal sheaf sequence \cite[p. 182]{Hartshorne} of
$X$ in $\PP^n$:
\begin{equation}\label{NormalSeq}
    0 \longrightarrow \Theta_X \longrightarrow i^{\ast}\Theta_{\PP^n} 
       \longrightarrow \mathcal{N}_{X/\PP^n} \longrightarrow 0,
\end{equation}
where $ \mathcal{N}_{X/\PP^n}\cong \bigoplus i^{\ast}\mathcal{O}_{\PP^n}(d_j) $.
Tensor the dual sequence of \eqref{NormalSeq} with the sheaf 
$\mathcal{L}_d^{\vee}= i^{\ast}\mathcal{O}_{\PP^n}(d-1)$ to obtain the exact sequence in $X$
$$
0 \longrightarrow \bigoplus i^{\ast} \mathcal{O}_{\PP^n}(-d_j+d-1)  \longrightarrow i^{\ast}\Omega^1_{\PP^n}(d-1) \longrightarrow \Omega^1_X \otimes i^{\ast} \mathcal{O}_{\PP^n}(d-1) \longrightarrow 0,
$$
whose long exact cohomology sequence 
\begin{equation}\label{GoalSeq}
    \begin{array}{c@{\hspace{3pt}}c@{\hspace{3pt}}c@{\hspace{3pt}}c@{\hspace{3pt}}c@{\hspace{3pt}}c@{\hspace{3pt}}c@{\hspace{3pt}}c}
         0 & \to & \bigoplus H^0(X,i^{\ast} \mathcal{O}_{\PP^n}(-d_j+d-1)) & \to &  H^0(X,i^{\ast}\Omega^1_{\PP^n}(d-1)) & \to & H^0(X,\Omega^1_X \otimes i^{\ast} \mathcal{O}_{\PP^n}(d-1)) & \to\\
         & \to & \bigoplus H^1(X,i^{\ast} \mathcal{O}_{\PP^n}(-d_j+d-1)) & \to & \cdots &  &  & 
    \end{array}
\end{equation}
involves the cohomology groups of interest.
Their dimensions  will be computed through those of
\textbf{(1)} 
 $h^k(X,  i^{\ast} \mathcal{O}_{\PP^n}(-d_j+d-1))$, for $ k=0,1 $, and 
\textbf{(2)}
 $ h^0(X,i^{\ast}\Omega^1_{\PP^n}(d-1)) $ in \eqref{GoalSeq}.
To this end, we use the following two results: 
\begin{proposition}\label{PropHart}
Let $X$ be a closed subset of $Y$, let $\mcE$ and $ \mathcal{F} $ be locally free sheaves on $X$ 
and $Y$ respectively,  and let $i : X \to Y$ be the inclusion. Then:
\begin{itemize}
    \item[(1) ] $H^k(X,\mcE)=H^k(Y,i_{\ast}\mcE)$, where $i_{\ast}\mcE$ is the extension of 
    $\mcE$ by zero outside $X$. In particular, 
    $ H^k(X,i^{\ast}\mathcal{F}) = H^k(Y,i_{\ast}i^{\ast}\mathcal{F}) $.
    \item[(2) ] The sequence 
    $$ 0 \longrightarrow \mathcal{F} \otimes_{\mathcal{O}_Y} \mathcal{I}_X 
     \longrightarrow \mathcal{F} \longrightarrow i_{\ast}i^{\ast} \mathcal{F} \longrightarrow 0 $$
    is exact.
\end{itemize}
\begin{proof}
For the proof of \textit{(1)}, see \cite[III Lemma 2.10]{Hartshorne}. 

Now we prove \textit{(2)}:
Consider the exact sequence on $Y$ of the ideal sheaf of $X$ \cite[p. 115]{Hartshorne}
$$
0 \longrightarrow \mathcal{I}_X \longrightarrow \mathcal{O}_Y \longrightarrow i_{\ast}\mathcal{O}_X \longrightarrow 0,
$$
and twist it by the sheaf $\mathcal{F} $ to obtain 
$$
0 \longrightarrow \mathcal{F} \otimes_{\mathcal{O}_Y} \mathcal{I}_X \longrightarrow \mathcal{F} 
    \longrightarrow \mathcal{F} \otimes_{\mathcal{O}_Y} i_{\ast}\mathcal{O}_X  \longrightarrow 0.
$$
It follows from the Projection Formula \cite[Ex. 5.1 (d)]{Hartshorne} that
$ \mathcal{F} \otimes_{\mathcal{O}_Y} i_{\ast}\mathcal{O}_X \cong i_{\ast}i^{\ast} \mathcal{F} $.
\end{proof}
\end{proposition}

We will also use the following computations of some cohomology groups in $ \PP^n $, which follow
from Bott's formula \cite[p. 4]{Okonek}:
\begin{lemma}\label{LemmBottFla}
\mbox{ } 
\begin{itemize}
    \item $h^q(\PP^n,\mathcal{O}_{\PP^n}(k))=0$, for $0<q<n$.
    \item $h^0(\PP^n,\mathcal{O}_{\PP^n}(k))=\binom{n+k}{k}$, for $k \geq 0$.
    \item $h^q(\PP^n, \Omega_{\PP^n}^1(k))=0$, for $0<q<n$, 
    except for $q=1$ and $k=0$, where it equals $1$.
    \item $h^0(\PP^n, \Omega_{\PP^n}^1(k))=(k-1)\binom{k+n-1}{k}$, for $k>1$.
\end{itemize}
\end{lemma} 
Now we begin with the computations: We start with the case of smooth quartic hypersurfaces in
$ \PP^3 $:
\begin{theorem}\label{Thm:Ex-1}
Let $ X \overset{ i }{\hookrightarrow} \PP^3 $
be a smooth quartic hypersurface,
and let $ d \geq 0 $ be an integer. Then
the dimensions of the spaces of sections that define foliations on $X$ with tangent sheaf 
$\mathcal{L}_d=i^{\ast}\mathcal{O}_{\PP^3}(1-d)$ are given by
\begin{equation*}
 h^0(X,\Omega^1_X \otimes i^{\ast} \mathcal{O}_{\PP^3}(d-1))
    = \begin{cases}
    0, & \text{if } d < 3,\\
    6, & \text{if } d = 3,\\
    20, & \text{if } d = 4,\\
    45, & \text{if }  d =5,\\
    4d^2-8d-16, & \text{if } d >5.
    \end{cases}
\end{equation*}
\begin{proof}
    From \eqref{NormalSeq}, 
    $\mathcal{N}_{X/\PP^3} \cong i^{\ast}\mathcal{O}_{\PP^3}(4)$ and the sequence \eqref{GoalSeq} becomes
    \begin{equation}\label{GoalSeq-1}
    \begin{array}{c@{\hspace{5pt}}c@{\hspace{5pt}}c@{\hspace{5pt}}c@{\hspace{5pt}}c@{\hspace{5pt}}c@{\hspace{5pt}}c@{\hspace{5pt}}c}
         0 & \to & H^0(X,i^{\ast} \mathcal{O}_{\PP^3}(d-5)) & \to &  H^0(X,i^{\ast}\Omega^1_{\PP^3}(d-1)) & \to & H^0(X,\Omega^1_X \otimes i^{\ast} \mathcal{O}_{\PP^3}(d-1)) & \to\\
         & \to & H^1(X,i^{\ast} \mathcal{O}_{\PP^3}(d-5)) & \to & \cdots &  &  & 
    \end{array}
\end{equation}
We compute the dimensions 
\textbf{(1)} $h^k(X, i^{\ast}\mathcal{O}_{\PP^3}(d-5)), \, k=0,1 ,$ 
and \textbf{(2)} $h^0(X, i^{\ast}\Omega_{\PP^3}^1(d-1))$ separately:

 \textbf{(1)}: By Corollary \ref{K3-CI-IdShf} (1) and 
 Proposition \ref{PropHart}, we have the short exact sequence
$$
    0 \longrightarrow \mathcal{O}_{\PP^3}(d-9) \longrightarrow \mathcal{O}_{\PP^3}(d-5) 
       \longrightarrow i_{\ast}i^{\ast}\mathcal{O}_{\PP^3}(d-5) \longrightarrow 0,
$$
whose long exact cohomology sequence is
\begin{equation*}\label{Sequence:Sequence1Example1}
    \begin{array}{c@{\hspace{5pt}}c@{\hspace{5pt}}c@{\hspace{5pt}}c@{\hspace{5pt}}c@{\hspace{5pt}}c@{\hspace{5pt}}c@{\hspace{5pt}}c}
         0 & \longrightarrow & H^0(\PP^3,\mathcal{O}_{\PP^3}(d-9)) & \longrightarrow &  H^0(\PP^3,\mathcal{O}_{\PP^3}(d-5)) & \longrightarrow & H^0(\PP^3,i_{\ast}i^{\ast}\mathcal{O}_{\PP^3}(d-5)) & \longrightarrow\\
         & \longrightarrow & H^1(\PP^3,\mathcal{O}_{\PP^3}(d-9)) & \longrightarrow &  H^1(\PP^3,\mathcal{O}_{\PP^3}(d-5)) & \longrightarrow & H^1(\PP^3,i_{\ast}i^{\ast}\mathcal{O}_{\PP^3}(d-5)) & \longrightarrow\\
         & \longrightarrow & H^2(\PP^3,\mathcal{O}_{\PP^3}(d-9)) & \longrightarrow &  \cdots& , & & 
    \end{array}
\end{equation*}
where 
$
 H^k(X, i^{\ast}\mathcal{O}_{\PP^3}(d-5)) \cong H^k(\PP^3, i_{\ast}i^{\ast}\mathcal{O}_{\PP^3}(d-5)), \, k=0,1.
$
By Lemma \ref{LemmBottFla}, we conclude that 
$ H^1(\PP^3,i_{\ast}i^{\ast}\mathcal{O}_{\PP^3}(d-5))=0 $ and moreover, that  
\begin{equation}\label{dim1Part1}
    h^0(X, i^{\ast}\mathcal{O}_{\PP^3}(d-5))= 
    \begin{cases}
    0, & \text{if } d < 5,\\
    1, & \text{if } d= 5,\\
    2d^2-20d+52, & \text{if } d >5.
    \end{cases}
\end{equation}
\textbf{(2)}: Now we come to $ h^0(X, i^{\ast}\Omega_{\PP^3}^1(d-1)) $:
Similar to case \textbf{(1)} above, we get the short exact sequence
$$
0 \longrightarrow \Omega_{\PP^3}^1(d-5) \longrightarrow \Omega_{\PP^3}^1(d-1) 
    \longrightarrow i_{\ast}i^{\ast}\Omega_{\PP^3}^1(d-1) \longrightarrow 0,
$$
whose long exact cohomology sequence is
\begin{equation}\label{SeqCotSeq-1}
    \begin{array}{c@{\hspace{5pt}}c@{\hspace{5pt}}c@{\hspace{5pt}}c@{\hspace{5pt}}c@{\hspace{5pt}}c@{\hspace{5pt}}c@{\hspace{5pt}}c}
         0 & \longrightarrow & H^0(\PP^3, \Omega_{\PP^3}^1(d-5)) & \longrightarrow &  H^0(\PP^3, \Omega_{\PP^3}^1(d-1)) & \longrightarrow & H^0(\PP^3, i_{\ast}i^{\ast}\Omega_{\PP^3}^1(d-1)) & \longrightarrow\\
         & \longrightarrow & H^1(\PP^3, \Omega_{\PP^3}^1(d-5) ) & \longrightarrow & H^1(\PP^3, \Omega_{\PP^3}^1(d-1) ) & \longrightarrow & \cdots &
    \end{array}
\end{equation}
where 
$
H^0(X, i^{\ast}\Omega_{\PP^3}^1(d-1)) \cong  H^0(\PP^3, i_{\ast}i^{\ast}\Omega_{\PP^3}^1(d-1)).
$
By Lemma \ref{LemmBottFla}, $H^1(\PP^3, \Omega_{\PP^3}^1(d-5))=0$ for all $d \neq 5$, 
and \eqref{SeqCotSeq-1} shows that
$$
h^0(\PP^3, i_{\ast}i^{\ast}\Omega_{\PP^3}^1(d-1)) = 
 h^0(\PP^3, \Omega_{\PP^3}^1(d-1))-h^0(\PP^3, \Omega_{\PP^3}^1(d-5)).
$$
For $ d=5 $, \eqref{SeqCotSeq-1} reduces to
$$
0 \longrightarrow H^0(\PP^3, \Omega_{\PP^3}^1(4)) \longrightarrow 
  H^0(\PP^3, i_{\ast}i^{\ast}\Omega_{\PP^3}^1(4)) \longrightarrow 
  H^1(\PP^3, \Omega_{\PP^3}^1) \longrightarrow 0,
$$
by Lemma \ref{LemmBottFla}.
The final result is
\begin{equation}\label{dim2Part1}
    h^0(X, i^{\ast}\Omega_{\PP^3}^1(d-1))= 
    \begin{cases}
    0, & \text{if } d < 3 ,\\
    6d^2-28d+36, & \text{if } d \geq 3.
    \end{cases}
\end{equation}
Since 
$ H^1(\PP^3,i_{\ast}i^{\ast}\mathcal{O}_{\PP^3}(d-5))= 0 $, the result follows from formulas \eqref{dim1Part1} and \eqref{dim2Part1} and the sequence \eqref{GoalSeq-1}.
\end{proof}
\end{theorem}
\begin{theorem}\label{Thm:Ex-2}
 Let $X \overset{ i }{\hookrightarrow} \PP^4$ be a smooth complete intersection surface of 
 type $(2,3)$ in $ \PP^4 $ and let $ d \geq 0 $ be an integer. Then 
 the dimensions of the spaces of sections that define foliations on $X$ with tangent sheaf 
 $\mathcal{L}_d=i^{\ast}\mathcal{O}_{\PP^4}(1-d)$ are given by
\begin{equation*}
    h^0(X, \Omega^1_X \otimes i^{\ast} \mathcal{O}_{\PP^4}(d-1)) = 
    \begin{cases}
    0, & \text{if } d < 3,\\
    10, & \text{if } d = 3,\\
    35, & \text{if } d = 4,\\
    6d^2-12d-14, & \text{if } d >4.
    \end{cases}
\end{equation*}
\begin{proof}
 From \eqref{NormalBundle}, 
 $\mathcal{N}_{X/\PP^4} \cong i^{\ast}(\mathcal{O}_{\PP^4}(2) \oplus \mathcal{O}_{\PP^4}(3))$, and the sequence \eqref{GoalSeq} becomes
    \begin{equation}\label{GoalSeq-2}
    \begin{array}{c@{\hspace{3pt}}c@{\hspace{3pt}}c@{\hspace{3pt}}c@{\hspace{3pt}}c@{\hspace{3pt}}c@{\hspace{3pt}}c@{\hspace{3pt}}c}
         0 & \to & H^0(X,i^{\ast} \mathcal{O}_{\PP^4}(d-3)) \oplus H^0(X,i^{\ast} \mathcal{O}_{\PP^4}(d-4)) & \to &  H^0(X,i^{\ast}\Omega^1_{\PP^4}(d-1)) & \to & H^0(X,\Omega^1_X \otimes i^{\ast} \mathcal{O}_{\PP^4}(d-1)) & \to\\
         & \to & H^1(X,i^{\ast} \mathcal{O}_{\PP^4}(d-3)) \oplus H^1(X,i^{\ast} \mathcal{O}_{\PP^4}(d-4)) & \to & \cdots &  &  & 
    \end{array}
\end{equation}
We will compute
\textbf{(1)}
$ h^k(X, i^{\ast}\mathcal{O}_{\PP^4}(d-3))$, $h^k(X, i^{\ast}\mathcal{O}_{\PP^4}(d-4)), \, 
k = 0,1 $, and 
\textbf{(2)} $ h^0(X, i^{\ast}\Omega_{\PP^4}^1(d-1))$, separately:

\textbf{(1)}: For
$ h^k(X, i^{\ast}\mathcal{O}_{\PP^4}(d-3)), \, k = 0,1 $:
Consider the sequence $(2)$ of Proposition \ref{PropHart}, with 
$ \mathcal{F} = \mathcal{O}_{\PP^4}(d-3) $. Apply Lemma \ref{aCM} to its 
associated cohomology sequence to obtain
\begin{equation}\label{Seq1-2}
 \begin{array}{c@{\hspace{5pt}}c@{\hspace{5pt}}c@{\hspace{5pt}}c@{\hspace{5pt}}c@{\hspace{5pt}} c@{\hspace{5pt}}c@{\hspace{5pt}}c}
    0 & \longrightarrow & H^0(\PP^4, \mathcal{I}_X(d-3)) & \longrightarrow &  
     H^0(\PP^4, \mathcal{O}_{\PP^4}(d-3)) & \longrightarrow & 
      H^0(\PP^4, i_{\ast}i^{\ast}\mathcal{O}_{\PP^4}(d-3)) & \longrightarrow\\
    & \longrightarrow & 0 & 
     \longrightarrow &  H^1(\PP^4, \mathcal{O}_{\PP^4}(d-3)) & \longrightarrow & 
      H^1(\PP^4, i_{\ast}i^{\ast}\mathcal{O}_{\PP^4}(d-3)) & \longrightarrow\\
        & \longrightarrow & 0 & \longrightarrow &  \cdots &  &  & 
    \end{array}
\end{equation}
where $H^k(\PP^4, i_{\ast}i^{\ast}\mathcal{O}_{\PP^4}(d-3)) 
 \cong H^k(X, i^{\ast}\mathcal{O}_{\PP^4}(d-3)), \, k = 0,1$.
Hence $ h^1(X, i^{\ast}\mathcal{O}_{\PP^4}(d-3)) = 0 $ by Lemma \ref{LemmBottFla}.

Now tensor the sequence $(2)$ of  Corollary \ref{K3-CI-IdShf}  with $\mathcal{O}_{\PP^4}(d-3)$
 to obtain
\begin{equation*}
    0 \longrightarrow \mathcal{O}_{\PP^4}(d-8) \longrightarrow \mathcal{O}_{\PP^4}(d-5)\oplus \mathcal{O}_{\PP^4}(d-6) \longrightarrow 
     \mathcal{I}_X(d-3) \longrightarrow 0.
\end{equation*}
By Lemma \ref{LemmBottFla}, a part of its long exact cohomology sequence is
\begin{equation*}
0 \to H^0(\PP^4,\mathcal{O}_{\PP^4}(d-8)) \to 
     H^0(\PP^4, \mathcal{O}_{\PP^4}(d-5)) \oplus H^0(\PP^4, \mathcal{O}_{\PP^4}(d-6)) 
     \to  H^0(\PP^4, \mathcal{I}_X(d-3)) \to 0.
\end{equation*}
This allows us to compute the value of 
$ h^0(\PP^4, \mathcal{I}_X(d-3)) $ and
the value of $ h^0(\PP^4, i_{\ast}i^{\ast}\mathcal{O}_{\PP^4}(d-3)) $, from \eqref{Seq1-2}.

In a completely analogous way one proves that 
$ h^1(\PP^4, i_{\ast}i^{\ast}\mathcal{O}_{\PP^4}(d-4))=0 $ and computes the
value of $ h^0(\PP^4, i_{\ast}i^{\ast}\mathcal{O}_{\PP^4}(d-4)) $. The final result is
\begin{equation}\label{dim1Part2}
    h^0(X,i^{\ast}\mathcal{O}_{\PP^4}(d-3)) + h^0(X,i^{\ast} \mathcal{O}_{\PP^4}(d-4))= 
    \begin{cases}
    0, & \text{if } d <3,\\
    1, & \text{if } d= 3,\\
    6, & \text{if } d= 4,\\
    6d^2-42d+79, & \text{if } d >4. 
    \end{cases}
\end{equation}
\textbf{(2)}: Now we come to
$h^0(X, i^{\ast}\Omega_{\PP^4}^1(d-1))$: we obtain from Proposition 
\ref{PropHart} the long exact cohomology sequence
\begin{equation}\label{SeqCotSeq-2}
 \begin{array}{c@{\hspace{5pt}}c@{\hspace{5pt}}c@{\hspace{5pt}}c@{\hspace{5pt}}c@{\hspace{5pt}}c@{\hspace{5pt}}c@{\hspace{5pt}}c@{\hspace{5pt}}c}
  0 & \longrightarrow & H^0(\PP^4, \Omega_{\PP^4}^1(d-1) \otimes \mathcal{I}_X ) & \longrightarrow &  H^0(\PP^4, \Omega_{\PP^4}^1(d-1)) 
  & \longrightarrow & H^0(\PP^4, i_{\ast}i^{\ast}\Omega_{\PP^4}^1(d-1)) & \longrightarrow &\\
  & \longrightarrow & H^1(\PP^4, \Omega_{\PP^4}^1(d-1) \otimes \mathcal{I}_X ) 
  & \longrightarrow & H^1(\PP^4, \Omega_{\PP^4}^1(d-1) ) & \longrightarrow 
  & H^1(\PP^4, i_{\ast}i^{\ast}\Omega_{\PP^4}^1(d-1)) & \longrightarrow & \cdots \\
    \end{array}
\end{equation}
where 
$H^0(\PP^4, i_{\ast}i^{\ast}\Omega_{\PP^4}^1(d-1)) 
\cong H^0(X,i^{\ast}\Omega_{\PP^4}^1(d-1))$. Tensor the sequence (2) of Corollary \ref{K3-CI-IdShf} with the sheaf 
$\Omega_{\PP^4}^1(d-1)$, to obtain
\begin{equation*}
    0 \longrightarrow \Omega_{\PP^4}^1(d-6) \longrightarrow \Omega_{\PP^4}^1(d-3) \oplus \Omega_{\PP^4}^1(d-4) \longrightarrow \Omega_{\PP^4}^1(d-1) \otimes \mathcal{I}_X \longrightarrow 0,
\end{equation*}
with cohomology sequence
\begin{equation*}
 \begin{array}{c@{\hspace{3pt}}c@{\hspace{3pt}}c@{\hspace{3pt}}c@{\hspace{3pt}}
   c@{\hspace{3pt}}c@{\hspace{3pt}}c@{\hspace{3pt}}c@{\hspace{3pt}}c}
    0 & \to & H^0(\PP^4, \Omega_{\PP^4}^1(d-6)) & \to 
      &  H^0(\PP^4, \Omega_{\PP^4}^1(d-3)) \oplus H^0(\PP^4, \Omega_{\PP^4}^1(d-4)) 
      & \to & H^0(\PP^4, \Omega_{\PP^4}^1(d-1) \otimes \mathcal{I}_X) & \to &\\
      & \to & H^1(\PP^4, \Omega_{\PP^4}^1(d-6)) & \to 
      &  H^1(\PP^4, \Omega_{\PP^4}^1(d-3)) \oplus H^1(\PP^4, \Omega_{\PP^4}^1(d-4)) 
      & \to & H^1(\PP^4, \Omega_{\PP^4}^1(d-1) \otimes \mathcal{I}_X) & \to & 0. 
    \end{array}
\end{equation*}
Applying Lemma \ref{LemmBottFla} to the previous cohomology sequences we obtain, 
for $d \neq 3,4,6 $, that
$$
h^0(\PP^4, i_{\ast}i^{\ast}\Omega_{\PP^4}^1(d-1)) = 
 h^0(\PP^4, \Omega_{\PP^4}^1(d-1))-h^0(\PP^4, \Omega_{\PP^4}^1(d-1) \otimes \mathcal{I}_X),
$$
and 
$$
h^0(\PP^4, \Omega_{\PP^4}^1(d-1) \otimes \mathcal{I}_X) = 
 h^0(\PP^4, \Omega_{\PP^4}^1(d-3)) + 
 h^0(\PP^4, \Omega_{\PP^4}^1(d-4))-h^0(\PP^4, \Omega_{\PP^4}^1(d-6)).
$$
When $d=3,4$ and $6$, Lemma \ref{LemmBottFla} shows that the sequence \eqref{SeqCotSeq-2} reduces respectively to
$$0 \longrightarrow H^0(\PP^4, \Omega_{\PP^4}^1(2)) \longrightarrow H^0(\PP^4, i_{\ast}i^{\ast}\Omega_{\PP^4}^1(2)) \longrightarrow H^1(\PP^4, \Omega_{\PP^4}^1(2)\otimes \mathcal{I}_X) \longrightarrow 0,$$
$$0 \longrightarrow H^0(\PP^4, \Omega_{\PP^4}^1(3)) \longrightarrow H^0(\PP^4, i_{\ast}i^{\ast}\Omega_{\PP^4}^1(3)) \longrightarrow H^1(\PP^4, \Omega_{\PP^4}^1(3)\otimes \mathcal{I}_X) \longrightarrow 0,$$
and
$$0 \longrightarrow H^0(\PP^4, \Omega_{\PP^4}^1(5)\otimes \mathcal{I}_X) \longrightarrow H^0(\PP^4, \Omega_{\PP^4}^1(5)) \longrightarrow H^0(\PP^4, i_{\ast}i^{\ast}\Omega_{\PP^4}^1(5))  \longrightarrow 0,$$
allowing us to compute the remaining dimensions. The final result is
\begin{equation}\label{dim2Part2}
h^0(X, i^{\ast}\Omega_{\PP^4}^1(d-1))= 
    \begin{cases}
    0, & \text{if } d < 3,\\
    12d^2-54d+65, & \text{if } d \geq 3.    
    \end{cases}
\end{equation}
Finally, since 
$ H^1(\PP^4, i_{\ast}i^{\ast}\mathcal{O}_{\PP^4}(d-3)) 
 \oplus H^1(\PP^4, i_{\ast}i^{\ast}\mathcal{O}_{\PP^4}(d-4))=0 $,
 the result follows from formulas
\eqref{dim1Part2} and \eqref{dim2Part2} and the sequence \eqref{GoalSeq-2}.
\end{proof}
\end{theorem}

The last case of projective complete intersection K3 surface is the following:
\begin{theorem}\label{Thm:Ex-3}
Let $X \overset{ i }{\hookrightarrow} \PP^5$ be a smooth complete intersection surface of type $(2,2,2)$ in $ \PP^5 $, and let $ d\geq 0$ be an integer. Then
the dimensions of the spaces of sections that define foliations on $X$ with tangent 
sheaf $\mathcal{L}_d=i^{\ast}\mathcal{O}_{\PP^5}(1-d)$ are given by
\begin{equation*}
    h^0(X, \Omega^1_X \otimes i^{\ast} \mathcal{O}_{\PP^5}(d-1)) = 
    \begin{cases}
    0, & \text{if } d < 3,\\
    15, & \text{if } d = 3,\\
    8d^2-16d-12, & \text{if } d >3.
    \end{cases}
\end{equation*}
\begin{proof}
  From \eqref{NormalBundle}, $\mathcal{N}_{X/\PP^4} \cong i^{\ast}\mathcal{O}_{\PP^5}(2)^{\oplus 3}$, 
 and the sequence \eqref{GoalSeq} becomes
    \begin{equation}\label{GoalSeq-3}
    \begin{array}{c@{\hspace{5pt}}c@{\hspace{5pt}}c@{\hspace{5pt}}c@{\hspace{5pt}}c@{\hspace{5pt}}c@{\hspace{5pt}}c@{\hspace{5pt}}c}
         0 & \longrightarrow & H^0(X,i^{\ast} \mathcal{O}_{\PP^5}(d-3))^{\oplus 3} & \longrightarrow &  H^0(X,i^{\ast}\Omega^1_{\PP^5}(d-1)) & \longrightarrow & H^0(X,\Omega^1_X \otimes i^{\ast} \mathcal{O}_{\PP^5}(d-1)) & \longrightarrow\\
         & \longrightarrow & H^1(X,i^{\ast} \mathcal{O}_{\PP^5}(d-3))^{\oplus 3} & \longrightarrow & \cdots &  &  & 
    \end{array}
\end{equation}
We compute the dimensions 
\textbf{(1)} $h^k(X, i^{\ast}\mathcal{O}_{\PP^5}(d-3)), k = 0, 1 $, and 
\textbf{(2)} $h^0(X, i^{\ast}\Omega_{\PP^5}^1(d-1))$ separately:

\textbf{(1)}: Consider the sequence $(2)$ of Proposition \ref{PropHart}, with 
$ \mathcal{F} = \mathcal{O}_{\PP^5}(d-3) $. Apply Lemma \ref{aCM} to its 
associated cohomology sequence to obtain 
\begin{equation}\label{Seq1-3}
\begin{array}{c@{\hspace{5pt}}c@{\hspace{5pt}}c@{\hspace{5pt}}c@{\hspace{5pt}}c@{\hspace{5pt}}c@{\hspace{5pt}}c@{\hspace{5pt}}c}
 0 & \longrightarrow & H^0(\PP^5,\mathcal{I}_X(d-3)) 
   & \longrightarrow &  H^0(\PP^5, \mathcal{O}_{\PP^5}(d-3)) 
   & \longrightarrow & H^0(\PP^5, i_{\ast}i^{\ast}\mathcal{O}_{\PP^5}(d-3)) 
   & \longrightarrow\\
    & \longrightarrow & 0
    & \longrightarrow &  H^1(\PP^5, \mathcal{O}_{\PP^5}(d-3)) 
    & \longrightarrow & H^1(\PP^5, i_{\ast}i^{\ast}\mathcal{O}_{\PP^5}(d-3)) 
    & \longrightarrow\\
     & \longrightarrow & 0 
     & \longrightarrow &  \cdots &  &  & 
    \end{array}
\end{equation}
where 
$ H^k(\PP^5, i_{\ast}i^{\ast}\mathcal{O}_{\PP^5}(d-3)) 
  \cong H^k(X, i^{\ast}\mathcal{O}_{\PP^5}(d-3)), k= 0,1 $.
Hence $ h^1(\PP^5, i^{\ast}\mathcal{O}_{\PP^5}(d-3))=0 $ by Lemma \ref{LemmBottFla}.  
  Now tensor the short exact sequence (3) of 
  Corollary \ref{K3-CI-IdShf} with $\mathcal{O}_{\PP^5}(d-3)$
  to obtain
\begin{equation*}
    0 \longrightarrow \left(\rfrac{\mathcal{O}_{\PP^5}(-4)^{\oplus 3}}{\mathcal{O}_{\PP^5}(-6)}\right) \otimes \mathcal{O}_{\PP^5}(d-3) \longrightarrow \mathcal{O}_{\PP^5}(d-5)^{\oplus 3} \longrightarrow \mathcal{O}_{\PP^5}(d-3) \otimes \mathcal{I}_X \longrightarrow 0.
\end{equation*}
After Lemma \ref{LemmBottFla}, its exact cohomology sequence is
\begin{equation}\label{SeqCohIdTen}
\begin{array}{c@{\hspace{2pt}}c@{\hspace{2pt}}c@{\hspace{2pt}}c@{\hspace{2pt}}c@{\hspace{2pt}}c@{\hspace{2pt}}c@{\hspace{2pt}}c}
 0 & \to 
 & H^0\left(\PP^5,\left(\rfrac{\mathcal{O}_{\PP^5}(-4)^{\oplus 3}}{\mathcal{O}_{\PP^5}(-6)}\right) \otimes \mathcal{O}_{\PP^5}(d-3)\right) 
 & \to &  H^0(\PP^5, \mathcal{O}_{\PP^5}(d-5))^{\oplus 3} & \to 
  & H^0(\PP^5, \mathcal{I}_X(d-3)) & \to
  \\
  & \to 
  & H^1\left(\PP^5,\left(\rfrac{\mathcal{O}_{\PP^5}(-4)^{\oplus 3}}{\mathcal{O}_{\PP^5}  (-6)}\right) \otimes \mathcal{O}_{\PP^5}(d-3)\right) 
  & \to & 0 & \to & \cdots
    \end{array}
\end{equation}
For the first column of \eqref{SeqCohIdTen}, tensor the short exact sequence
\begin{equation}\label{New}
    0 \longrightarrow \mathcal{O}_{\PP^5}(-6) \longrightarrow \mathcal{O}_{\PP^5}(-4)^{\oplus 3} \longrightarrow \rfrac{\mathcal{O}_{\PP^5}(-4)^{\oplus 3}}{\mathcal{O}_{\PP^5}(-6)} \longrightarrow 0,
\end{equation}
with $\mathcal{O}_{\PP^5}(d-3)$. After Lemma \ref{LemmBottFla}, a part of
its exact cohomology sequence is 
\begin{equation}\label{NewSeq:Cohom}
\begin{array}{c@{\hspace{1pt}}c@{\hspace{1pt}}c@{\hspace{1pt}}c@{\hspace{1pt}}c@{\hspace{1pt}}c@{\hspace{1pt}}c@{\hspace{1pt}}c}
0 & \to & H^0(\PP^5,\mathcal{O}_{\PP^5}(d-9)) & \to 
  & H^0(\PP^5, \mathcal{O}_{\PP^5}(d-7))^{\oplus 3} & \to 
  & H^0\left(\PP^5,\left(\rfrac{\mathcal{O}_{\PP^5}(-4)^{\oplus 3}}{\mathcal{O}_{\PP^5}(-6)}\right) \otimes \mathcal{O}_{\PP^5}(d-3)\right)  & \to
  \\
& \to & 0 & \to
& 0 & \to 
   & H^1\left(\PP^5,\left(\rfrac{\mathcal{O}_{\PP^5}(-4)^{\oplus 3}}{\mathcal{O}_{\PP^5}(-6)}\right) \otimes \mathcal{O}_{\PP^5}(d-3)\right) &  \to
   \\
& \to & 0 & \to
& \cdots
    \end{array}
\end{equation}

 In sum, we can now compute from \eqref{Seq1-3} the value of 
 $ h^0(\PP^5, i_{\ast}i^{\ast}\mathcal{O}_{\PP^5}(d-3)) $
 using the first rows of \eqref{NewSeq:Cohom} and \eqref{SeqCohIdTen}.
 The final result is

\begin{equation}\label{dim1Part3}
    3 \cdot h^0(X,i^{\ast}\mathcal{O}_{\PP^5}(d-3)) = 
    \begin{cases}
    0, & \text{if } d <3,\\
    3, & \text{if } d= 3,\\
    12d^2-72d+114, & \text{if } d >3.  
    \end{cases}
\end{equation}
\textbf{(2)}: Now we come to the computation of $h^0(X, i^{\ast}\Omega_{\PP^5}^1(d-1))$: From 
Proposition \ref{PropHart} we obtain the long exact cohomology sequence
\begin{equation}\label{SeqCotSeq-3}
    \begin{array}{c@{\hspace{5pt}}c@{\hspace{5pt}}c@{\hspace{5pt}}c@{\hspace{5pt}}c@{\hspace{5pt}}c@{\hspace{5pt}}c@{\hspace{5pt}}c}
         0 & \longrightarrow & H^0(\PP^5, \Omega_{\PP^5}^1(d-1) \otimes \mathcal{I}_X ) & \longrightarrow &  H^0(\PP^5, \Omega_{\PP^5}^1(d-1)) & \longrightarrow & H^0(\PP^5, i_{\ast}i^{\ast}\Omega_{\PP^5}^1(d-1)) & \longrightarrow\\
         & \longrightarrow & H^1(\PP^5, \Omega_{\PP^5}^1(d-1) \otimes \mathcal{I}_X ) & \longrightarrow & H^1(\PP^5, \Omega_{\PP^5}^1(d-1) ) & \longrightarrow & H^1(\PP^5, i_{\ast}i^{\ast}\Omega_{\PP^5}^1(d-1)) & \longrightarrow \\
         & \longrightarrow & H^2(\PP^5, \Omega_{\PP^5}^1(d-1) \otimes \mathcal{I}_X ) & \longrightarrow & H^2(\PP^5, \Omega_{\PP^5}^1(d-1)) & \longrightarrow & \cdots &
    \end{array}
\end{equation}
where $H^0(\PP^5, i_{\ast}i^{\ast}\Omega_{\PP^5}^1(d-1)) 
 \cong H^0(X, i^{\ast}\Omega_{\PP^5}^1(d-1))$. 
 Tensor the short exact sequence (3) from Corollary 
\ref{K3-CI-IdShf} with $\Omega_{\PP^5}^1(d-1)$ to obtain
\begin{equation*}
    0 \longrightarrow \left(\rfrac{\mathcal{O}_{\PP^5}(-4)^{\oplus 3}}{\mathcal{O}_{\PP^5}(-6)}\right) \otimes \Omega_{\PP^5}^1(d-1) \longrightarrow \Omega_{\PP^5}^1(d-3)^{\oplus 3} \longrightarrow \Omega_{\PP^5}^1(d-1) \otimes \mathcal{I}_X \longrightarrow 0,
\end{equation*}
whose long exact cohomology sequence is
\begin{equation}\label{SeqCohIdTenP2}
 \begin{array}{c@{\hspace{3pt}}c@{\hspace{3pt}}c@{\hspace{3pt}}c@{\hspace{3pt}}c@{\hspace{3pt}}
  c@{\hspace{3pt}}c@{\hspace{3pt}}c}
   0 & \to 
     & H^0\left(\PP^5,\left(\rfrac{\mathcal{O}_{\PP^5}(-4)^{\oplus 3}}{\mathcal{O}_{\PP^5}(-6)}\right) \otimes \Omega_{\PP^5}^1(d-1) \right) 
     & \to &  H^0(\PP^5, \Omega_{\PP^5}^1(d-3))^{\oplus 3} 
     & \to & H^0(\PP^5, \Omega_{\PP^5}^1(d-1) \otimes \mathcal{I}_X)) & \to\\
     & \to 
     & H^1\left(\PP^5,\left(\rfrac{\mathcal{O}_{\PP^5}(-4)^{\oplus 3}}{\mathcal{O}_{\PP^5}(-6)}\right) \otimes \Omega_{\PP^5}^1(d-1) \right) 
     & \to &  H^1(\PP^5, \Omega_{\PP^5}^1(d-3))^{\oplus 3} 
     & \to & H^1(\PP^5, \Omega_{\PP^5}^1(d-1) \otimes \mathcal{I}_X)) & \to\\
     & \to 
     & H^2\left(\PP^5,\left(\rfrac{\mathcal{O}_{\PP^5}(-4)^{\oplus 3}}{\mathcal{O}_{\PP^5}(-6)}\right) \otimes \Omega_{\PP^5}^1(d-1) \right) 
     & \to &  H^2(\PP^5, \Omega_{\PP^5}^1(d-3))^{\oplus 3} 
     & \to & H^2(\PP^5, \Omega_{\PP^5}^1(d-1) \otimes \mathcal{I}_X)) & \to\\
     & \to & H^3\left(\PP^5,\left(\rfrac{\mathcal{O}_{\PP^5}(-4)^{\oplus 3}}{\mathcal{O}_{\PP^5}(-6)}\right) \otimes \Omega_{\PP^5}^1(d-1) \right) & \to &  \cdots &  &  & 
    \end{array}
\end{equation}
For the first column of \eqref{SeqCohIdTenP2}, tensor the sequence \eqref{New}
with the sheaf $\Omega_{\PP^5}^1(d-1)$, and consider its long exact cohomology sequence:
\begin{equation*}
    \begin{array}{c@{\hspace{5pt}}c@{\hspace{5pt}}c@{\hspace{5pt}}c@{\hspace{5pt}}c@{\hspace{5pt}}c@{\hspace{5pt}}c@{\hspace{5pt}}c}
         0 & \to & H^0(\PP^5,\Omega_{\PP^5}^1(d-7)) & \to &  H^0(\PP^5, \Omega_{\PP^5}^1(d-5))^{\oplus 3} & \to & H^0\left(\PP^5,\left(\rfrac{\mathcal{O}_{\PP^5}(-4)^{\oplus 3}}{\mathcal{O}_{\PP^5}(-6)}\right) \otimes \Omega_{\PP^5}^1(d-1)\right) & \to\\
         & \to & H^1(\PP^5,\Omega_{\PP^5}^1(d-7)) & \to &  H^1(\PP^5, \Omega_{\PP^5}^1(d-5))^{\oplus 3} & \to & H^1\left(\PP^5,\left(\rfrac{\mathcal{O}_{\PP^5}(-4)^{\oplus 3}}{\mathcal{O}_{\PP^5}(-6)}\right) \otimes \Omega_{\PP^5}^1(d-1)\right) & \to\\
         & \to & H^2(\PP^5,\Omega_{\PP^5}^1(d-7)) & \to &  H^2(\PP^5, \Omega_{\PP^5}^1(d-5))^{\oplus 3} & \to & H^2\left(\PP^5,\left(\rfrac{\mathcal{O}_{\PP^5}(-4)^{\oplus 3}}{\mathcal{O}_{\PP^5}(-6)}\right) \otimes \Omega_{\PP^5}^1(d-1)\right) & \to\\
         & \to & H^3(\PP^5,\Omega_{\PP^5}^1(d-7)) & \to &  H^3(\PP^5, \Omega_{\PP^5}^1(d-5))^{\oplus 3} & \to & H^3\left(\PP^5,\left(\rfrac{\mathcal{O}_{\PP^5}(-4)^{\oplus 3}}{\mathcal{O}_{\PP^5}(-6)}\right) \otimes \Omega_{\PP^5}^1(d-1)\right) & \to\\
         & \to & H^4(\PP^5,\Omega_{\PP^5}^1(d-7)) & \to &  \cdots &  & & 
    \end{array}
\end{equation*}

Applying Lemma \ref{LemmBottFla} to the previous sequences we obtain, for $ d \neq 3,5,7 $,
that
$$h^0(\PP^5, i_{\ast}i^{\ast}\Omega_{\PP^5}^1(d-1))=h^0(\PP^5, \Omega_{\PP^5}^1(d-1))-h^0(\PP^5, \Omega_{\PP^5}^1(d-1) \otimes \mathcal{I}_X),$$
and 
$$h^0(\PP^5, \Omega_{\PP^5}^1(d-1) \otimes \mathcal{I}_X)=3h^0(\PP^5, \Omega_{\PP^5}^1(d-3))-(3h^0(\PP^5, \Omega_{\PP^5}^1(d-5))-h^0(\PP^5, \Omega_{\PP^5}^1(d-7))).$$

When $d=3,5$ and $7$, Lemma \ref{LemmBottFla} shows that the sequence \eqref{SeqCotSeq-3} reduces respectively to
$$0 \longrightarrow H^0(\PP^5, \Omega_{\PP^5}^1(2)) \longrightarrow H^0(\PP^5, i_{\ast}i^{\ast}\Omega_{\PP^5}^1(2)) \longrightarrow H^1(\PP^5, \Omega_{\PP^5}^1(2)\otimes \mathcal{I}_X) \longrightarrow 0,$$
$$0 \longrightarrow H^0(\PP^5, \Omega_{\PP^5}^1(4)\otimes \mathcal{I}_X) \longrightarrow H^0(\PP^5, \Omega_{\PP^5}^1(4)) \longrightarrow H^0(\PP^5, i_{\ast}i^{\ast}\Omega_{\PP^5}^1(4))  \longrightarrow 0,$$
and
$$0 \longrightarrow H^0(\PP^5, \Omega_{\PP^5}^1(6)\otimes \mathcal{I}_X) \longrightarrow H^0(\PP^5, \Omega_{\PP^5}^1(6)) \longrightarrow H^0(\PP^5, i_{\ast}i^{\ast}\Omega_{\PP^5}^1(6))  \longrightarrow 0,$$
allowing us to compute the remaining dimensions. The final result is
\begin{equation}\label{dim2Part3}
h^0(X, i^{\ast}\Omega_{\PP^5}^1(d-1))= 
    \begin{cases}
    0, & \text{if } d < 3,\\
    20d^2-88d+102, & \text{if } d \geq 3.    
    \end{cases}
\end{equation}
Finally, since $H^1(\PP^5, i_{\ast}i^{\ast}\mathcal{O}_{\PP^5}(d-3)) =0$, the result follows
from formulas \eqref{dim1Part3} and \eqref{dim2Part3} and the sequence \eqref{GoalSeq-3}.
\end{proof}
\end{theorem}
\section{Proof of Theorem B}\label{Ch:PffThmB}
\textbf{Main idea for the proof of Theorem B:}
Let $ X $ be a K3 surface, $ \mathcal{L} $ an invertible sheaf on $ X $
and let $ s\in H^0(X, \Theta_X \otimes \mathcal{L}^{\vee}) $ 
be a 
section with isolated singularities, with singular scheme $ Z = Z_s $. 
The Koszul complex associated to $s$ \cite[p. 76]{RRAlgebra} gives the following resolution 
of the ideal sheaf $\mathcal{I}_Z \subset \mcO_X$ of $Z$
\begin{equation*}
    0 \longrightarrow \bigwedge^2 \left(\Omega^1_X \otimes \mathcal{L}\right)  \overset{i_s}{\longrightarrow} \bigwedge^1 \left(\Omega^1_X \otimes \mathcal{L}\right) \overset{i_s}{\longrightarrow} \mathcal{I}_Z \longrightarrow 0, 
\end{equation*}
where the morphisms $\iota_s$ are \textit{evaluation in} (or \textit{contraction by})
the section $s$. This sequence is equal to  
\begin{equation*}
    0 \longrightarrow \mathcal{L}^{\otimes 2} \overset{i_s}{\longrightarrow} \Omega^1_X \otimes \mathcal{L} \overset{i_s}{\longrightarrow} \mathcal{I}_Z \longrightarrow 0.
\end{equation*}
 Now tensor the sequence above with the sheaf 
 $\mcE=\Theta_X \otimes \mathcal{L}^{\vee}$, to obtain the exact sequence
 \begin{equation}\label{basicseq}
    0 \longrightarrow \Theta_X \otimes \mathcal{L} 
       \longrightarrow \mathscr{H}om(\Theta_X,\Theta_X)) 
        \overset{i_s \otimes \mathrm{I}_{\mcE}}{\longrightarrow} \mathcal{I}_Z \otimes \Theta_X \otimes \mathcal{L}^{\vee}       \longrightarrow 0.
\end{equation}
  Now recall that $\Theta_X$ is \textit{simple} \cite[Proposition 4.5]{Huybrechts}:
 $ H^0(X, \mathscr{H}om(\Theta_X,\Theta_X)) \cong \CC \cdot \mathrm{I}_{\Theta_X} $; 
 recall moreover that
 $ H^0(X, \mathcal{I}_Z \otimes \Theta_X \otimes \mathcal{L}^{\vee}) $ is the
 space of sections in $ H^0(X,\Theta_X \otimes \mathcal{L}^{\vee}) $
 that vanish in $ Z $ and finally, that
 the sheaf mapping $ \iota_s \otimes \mathrm{I}_{\mcE} $ is $ e \mapsto e( s ) $.

 Then, the long exact cohomology sequence associated to \eqref{basicseq} is
\begin{equation}\label{KoszCoh}
    \begin{array}{c@{\hspace{5pt}}c@{\hspace{5pt}}c@{\hspace{5pt}}c@{\hspace{5pt}}c@{\hspace{5pt}}c@{\hspace{5pt}}c@{\hspace{5pt}}c}
         0 & \longrightarrow & 0 & \longrightarrow &  \CC & \overset{\Phi}{\longrightarrow} & H^0(X, \mathcal{I}_Z \otimes \Theta_X \otimes \mathcal{L}^{\vee}) & \longrightarrow\\
         & \longrightarrow & H^1(X, \Theta_X \otimes \mathcal{L}) & \longrightarrow & H^1(X, \mathscr{H}om(\Theta_X,\Theta_X)) & \longrightarrow & \cdots & 
    \end{array}
\end{equation}
because the linear map $ \Phi = (i_s \otimes \mathrm{I}_{\mathscr{E}})^0 $
induced by $ \iota_s \otimes \mathrm{I}_{\mathscr{E}} $ in global sections
is the (obviously injective) map $ \Phi( \lambda\cdot\mathrm{I}_{\Theta_X} ) = \lambda\cdot s $.

It follows that the unique sections 
$ s' \in H^0(X, \mathcal{I}_Z \otimes \Theta_X \otimes \mathcal{L}^{\vee}) $
that vanish in $ Z = Z_s $ are those of the form $ s' = \lambda\cdot s $ if and only if
the map $ \Phi $ is surjective: this is, that $ \mcF = [ s ] $ is the unique foliation that
vanishes in $ Z = Z_s $ if and only if $ \Phi $ is surjective. From \eqref{KoszCoh},
this last condition holds if $H^1(X, \Theta_X \otimes \mathcal{L}) = 0 $.

In sum, \textbf{
if $H^1(X, \Theta_X \otimes \mathcal{L}) = 0$ then 
a foliation $ \mcF = [ s ] $ in $ X $
with isolated singularities and tangent sheaf $ \mcL $ is uniquely determined by its singular 
scheme $ Z = Z_s $. }

The statement in Theorem B is the result of the computation of those $d \in \ZZ $ for which
\begin{equation}\label{goal}
   H^1(X, \Theta_X \otimes i^{\ast}\mathcal{O}_{\PP^n}(1-d)) = 0, 
\end{equation}
for the three types of projective complete intersection K3 surfaces $ X $ 
(Lemma \ref{K3-CI}).

As in Chapter \ref{Ch:PffThmA}, the computations will follow a common pattern:
Tensor the normal sheaf sequence \eqref{NormalSeq}, 
with the sheaf $\mathcal{L}_d= i^{\ast}\mathcal{O}_{\PP^n}(1-d)$
to obtain the sequence on $X$:
$$
0 \longrightarrow \Theta_X \otimes i^{\ast}\mathcal{O}_{\PP^n}(1-d) 
   \longrightarrow i^{\ast}\Theta_{\PP^n}(1-d) 
    \longrightarrow \bigoplus i^{\ast}\mathcal{O}_{\PP^n}(d_j-d+1) 
     \longrightarrow 0,
$$
\medskip
whose long exact cohomology sequence is
\begin{equation}\label{NorCohSeqV2}
    \begin{array}{c@{\hspace{5pt}}c@{\hspace{5pt}}c@{\hspace{5pt}}c@{\hspace{5pt}}c@{\hspace{5pt}}c@{\hspace{5pt}}c@{\hspace{5pt}}c}
         0 & \to & H^0(X, \Theta_X \otimes i^{\ast}\mathcal{O}_{\PP^n}(1-d)) & \to 
            &  H^0(X, i^{\ast}\Theta_{\PP^n}(1-d)) & \to 
             & \bigoplus H^0(X,i^{\ast}\mathcal{O}_{\PP^n}(d_j-d+1)) & \to \\
              & \to & H^1(X, \Theta_X \otimes i^{\ast}\mathcal{O}_{\PP^n}(1-d)) & \to 
               & H^1(X, i^{\ast}\Theta_{\PP^n}(1-d)) & \to & \cdots
    \end{array}
\end{equation}
The values of $ d $ for which both  
$h^0(X,i^{\ast}\mathcal{O}_{\PP^n}(d_j-d+1)) = h^1(X, i^{\ast}\Theta_{\PP^n}(1-d)) = 0 $,
for every $ d_j $, give our solution to \eqref{goal}.

As in the Introduction, the space of foliations 
$ \PP H^0( X, \mathscr{H}om(i^{\ast}\mathcal{O}_{\PP^n}(1-d), \Theta_X )) $
on $ X $ with tangent sheaf $\mathcal{L}_d=i^{\ast}\mathcal{O}_{\PP^n}(1-d)$ 
$ i^{\ast}\mathcal{O}_{\PP^n}(1-d) $
will be denoted by $ \mathscr{F}ol( i^{\ast}\mathcal{O}_{\PP^n}(1-d), X ) $.
\begin{theorem}\label{Thm:Uniq-1}
  Let $ X \overset{ i }{\hookrightarrow} \PP^3 $
  be a smooth quartic hypersurface, 
  and let $ d \geq 3 $ be an integer. Then for $d \geq 6$,
  any foliation $ \mathscr{F} = [s] \in \mathscr{F}ol( i^{\ast}\mathcal{O}_{\PP^3}(1-d), X ) $  
  with isolated singularities is uniquely determined by its singular scheme. 

\begin{proof}
    Since $\mathcal{N}_{X/\PP^3} \cong i^{\ast}\mathcal{O}_{\PP^3}(4)$, the sequence \eqref{NorCohSeqV2} becomes
    \begin{equation*}
    \begin{array}{c@{\hspace{5pt}}c@{\hspace{5pt}}c@{\hspace{5pt}}c@{\hspace{5pt}}c@{\hspace{5pt}}c@{\hspace{5pt}}c@{\hspace{5pt}}c}
         0 & \to & H^0(X, \Theta_X \otimes i^{\ast}\mathcal{O}_{\PP^3}(1-d)) & \to &  H^0(X, i^{\ast}\Theta_{\PP^3}(1-d)) & \to & H^0(X,i^{\ast}\mathcal{O}_{\PP^3}(5-d)) & \to \\
         & \to & H^1(X, \Theta_X \otimes i^{\ast}\mathcal{O}_{\PP^3}(1-d)) & \to & H^1(X, i^{\ast}\Theta_{\PP^3}(1-d)) & \to & \cdots
    \end{array}
\end{equation*}

To compute $h^0(X, i^{\ast}\mathcal{O}_{\PP^3}(5-d))$, observe that 
Proposition \ref{PropHart} yields the short exact sequence in $\PP^3$ 
    $$0 \longrightarrow \mathcal{O}_{\PP^3}(1-d) \longrightarrow \mathcal{O}_{\PP^3}(5-d) \longrightarrow i_{\ast}i^{\ast}\mathcal{O}_{\PP^3}(5-d) \longrightarrow 0,$$
which induces the long exact sequence
\begin{equation}\label{Seq2-1}
    \begin{array}{c@{\hspace{5pt}}c@{\hspace{5pt}}c@{\hspace{5pt}}c@{\hspace{5pt}}c@{\hspace{5pt}}c@{\hspace{5pt}}c@{\hspace{5pt}}c}
         0 & \longrightarrow & H^0(\PP^3,\mathcal{O}_{\PP^3}(1-d)) & \longrightarrow &  H^0(\PP^3,\mathcal{O}_{\PP^3}(5-d)) & \longrightarrow & H^0(\PP^3,i_{\ast}i^{\ast}\mathcal{O}_{\PP^3}(5-d)) & \longrightarrow\\
         & \longrightarrow & H^1(\PP^3,\mathcal{O}_{\PP^3}(1-d)) & \longrightarrow &  \cdots & &  & 
    \end{array}
\end{equation}
where $H^0(X, i^{\ast}\mathcal{O}_{\PP^3}(5-d)) \cong 
 H^0(\PP^3, i_{\ast}i^{\ast}\mathcal{O}_{\PP^3}(5-d))$. 
Since $H^k(\PP^3,\mathcal{O}_{\PP^3}(1-d))=0$ for $ k = 0, 1 $
by Lemma \ref{LemmBottFla}, we obtain 
\begin{equation}\label{dim1ex1}
    h^0(X, i^{\ast}\mathcal{O}_{\PP^3}(5-d))= 
    \begin{cases}
    2d^2-20d+52, & \text{if } d = 3, 4,\\
    1, & \text{if } d=5,\\
    0, & \text{if } d \geq 6.
    \end{cases}
\end{equation}
 Now we come to $h^1(X, i^{\ast}\Theta_{\PP^3}(1-d))$: consider the short exact sequence in 
 $\PP^3$
$$
0 \longrightarrow \Theta_{\PP^3}(-3-d) \longrightarrow \Theta_{\PP^3}(1-d) 
  \longrightarrow i_{\ast}i^{\ast}\Theta_{\PP^3}(1-d) \longrightarrow 0,
$$
given by Proposition \ref{PropHart}. This sequence induces the long exact sequence
\begin{equation}\label{TangSeq-1}
    \begin{array}{c@{\hspace{5pt}}c@{\hspace{5pt}}c@{\hspace{5pt}}c@{\hspace{5pt}}c@{\hspace{5pt}}c@{\hspace{5pt}}c@{\hspace{5pt}}c}
         0 & \longrightarrow & H^0(\PP^3, \Theta_{\PP^3}(-3-d)) & \longrightarrow 
         &  H^0(\PP^3, \Theta_{\PP^3}(1-d)) & \longrightarrow 
         & H^0(\PP^3, i_{\ast}i^{\ast}\Theta_{\PP^3}(1-d)) & \longrightarrow\\
         & \longrightarrow & H^1(\PP^3, \Theta_{\PP^3}(-3-d)) & \longrightarrow 
         &  H^1(\PP^3, \Theta_{\PP^3}(1-d)) & \longrightarrow 
         & H^1(\PP^3, i_{\ast}i^{\ast}\Theta_{\PP^3}(1-d)) & \longrightarrow\\
         & \longrightarrow & H^2(\PP^3, \Theta_{\PP^3}(-3-d)) & \longrightarrow & \cdots &  &  &
    \end{array}
\end{equation}
where $H^1(X, i^{\ast}\Theta_{\PP^3}(1-d)) \cong  
H^1(\PP^3, i_{\ast}i^{\ast}\Theta_{\PP^3}(1-d))$. Applying Serre duality and Lemma \ref{LemmBottFla}, we find that 
$H^1(\PP^3, \Theta_{\PP^3}(1-d))=0$, and that
$H^2(\PP^3, \Theta_{\PP^3}(-3-d))=0$ if and only if $d \neq 1$. Hence

\begin{equation}\label{dneq1}
    H^1(X,i^{\ast}\Theta_{\PP^3}(1-d)) = 0 \textup{ if and only if } d \neq 1,
\end{equation}

Combine \eqref{dim1ex1} and \eqref{dneq1} to conclude that 
$H^1(X, \Theta_X \otimes i^{\ast}\mathcal{O}_{\PP^3}(1-d))=0$ for $d \geq 6$, 
which finishes the proof. 
\end{proof}
\end{theorem}

\begin{theorem}\label{Thm:Uniq-2}
    Let $ X \overset{ i }{\hookrightarrow} \PP^4 $
    be a smooth complete intersection surface of type $(2,3)$ and let $ d\geq 3 $ be an integer.
    Then, for $ d \geq 5 $, 
    any foliation $ \mathscr{F} = [s] \in \mathscr{F}ol( i^{\ast}\mathcal{O}_{\PP^4}(1-d), X ) $  
    with isolated singularities is uniquely determined by its singular scheme.     

\begin{proof}
 Since $\mathcal{N}_{X/\PP^4} \cong i^{\ast}(\mathcal{O}_{\PP^4}(2) \oplus \mathcal{O}_{\PP^4}(3))$, the sequence \eqref{NorCohSeqV2} becomes
    \begin{equation*}
    \begin{array}{c@{\hspace{3pt}}c@{\hspace{3pt}}c@{\hspace{3pt}}c@{\hspace{3pt}}c@{\hspace{3pt}}c@{\hspace{3pt}}c@{\hspace{3pt}}c}
     0 & \to & H^0(X, \Theta_X \otimes i^{\ast}\mathcal{O}_{\PP^4}(1-d)) & \to &  H^0(X, i^{\ast}\Theta_{\PP^4}(1-d)) & \to & H^0(X,i^{\ast}\mathcal{O}_{\PP^4}(3-d)) \oplus H^0(X,i^{\ast}\mathcal{O}_{\PP^4}(4-d)) & \to \\
         & \to & H^1(X, \Theta_X \otimes i^{\ast}\mathcal{O}_{\PP^4}(1-d)) & \to &  H^1(X, i^{\ast}\Theta_{\PP^4}(1-d)) & \to & \cdots
    \end{array}
\end{equation*}
We compute the dimensions \textbf{(1)}
$h^0(X, i^{\ast}\mathcal{O}_{\PP^4}(3-d))$, $h^0(X, i^{\ast}\mathcal{O}_{\PP^4}(4-d))$  
and \textbf{(2)} $h^1(X, i^{\ast}\Theta_{\PP^4}(1-d))$, separately:

\textbf{(1)}:
For $h^0(X, i^{\ast}\mathcal{O}_{\PP^4}(3-d))$:
Consider the sequence $(2)$ of Proposition \ref{PropHart}, with 
$ \mathcal{F} = \mathcal{O}_{\PP^4}(3-d) $. Apply Lemma \ref{aCM} to its 
associated cohomology sequence to obtain

\begin{equation}\label{Seq2-2}
0  \longrightarrow  H^0(\PP^4,\mathcal{I}_X(3-d)) 
   \longrightarrow   H^0(\PP^4, \mathcal{O}_{\PP^4}(3-d)) 
   \longrightarrow  H^0(\PP^4, i_{\ast}i^{\ast}\mathcal{O}_{\PP^4}(3-d)) 
   \longrightarrow  0 
\end{equation}
where $H^0(\PP^4, i_{\ast}i^{\ast}\mathcal{O}_{\PP^4}(3-d)) 
\cong H^0(X, i^{\ast}\mathcal{O}_{\PP^4}(3-d))$. 

 Now, it is obvious that $h^0(\PP^4,\mathcal{I}_X(3-d)) = 0 $ but anyway we
 prove it:
tensor the sequence (2) 
in Corollary \ref{K3-CI-IdShf} with the sheaf $\mathcal{O}_{\PP^4}(3-d)$ to obtain
\begin{equation*}
0 \longrightarrow \mathcal{O}_{\PP^4}(-2-d) \longrightarrow 
   \mathcal{O}_{\PP^4}(1-d)\oplus \mathcal{O}_{\PP^4}(-d) 
   \longrightarrow \mathcal{I}_X(3-d) \longrightarrow 0.
\end{equation*}
The conclusion follows by applying Lemma \ref{LemmBottFla} to its 
associated cohomology sequence.

In an analogous way, one computes the
value of $h^0(X, i^{\ast}\mathcal{O}_{\PP^4}(4-d))$. The final result is
\begin{equation}\label{dim1ex2}
    h^0(X, i^{\ast}\mathcal{O}_{\PP^4}(3-d))+h^0(X, i^{\ast}\mathcal{O}_{\PP^4}(4-d))= 
    \begin{cases}
    6, & \text{if } d=3,\\
    1, & \text{if } d=4,\\
    0, & \text{if } d \geq 5.
    \end{cases}
\end{equation}
\textbf{(2)}: Now we come to
$h^1(X, i^{\ast}\Theta_{\PP^4}(1-d))$: Consider the short exact sequence in $\PP^4$
$$0 \longrightarrow \Theta_{\PP^4}(1-d) \otimes \mathcal{I}_X \longrightarrow \Theta_{\PP^4}(1-d) \longrightarrow i_{\ast}i^{\ast}\Theta_{\PP^4}(1-d) \longrightarrow 0,$$
given by Proposition \ref{PropHart}. This sequence induces the long exact sequence
\begin{equation}\label{TangSeq-2}
    \begin{array}{c@{\hspace{5pt}}c@{\hspace{5pt}}c@{\hspace{5pt}}c@{\hspace{5pt}}c@{\hspace{5pt}}c@{\hspace{5pt}}c@{\hspace{5pt}}c}
         0 & \longrightarrow & H^0(\PP^4, \Theta_{\PP^4}(1-d) \otimes \mathcal{I}_X) & \longrightarrow &  H^0(\PP^4, \Theta_{\PP^4}(1-d)) & \longrightarrow 
         & H^0(\PP^4, i_{\ast}i^{\ast}\Theta_{\PP^4}(1-d)) & \longrightarrow\\
         & \longrightarrow & H^1(\PP^4, \Theta_{\PP^4}(1-d) \otimes \mathcal{I}_X) 
         & \longrightarrow &  H^1(\PP^4, \Theta_{\PP^4}(1-d)) 
         & \longrightarrow & H^1(\PP^4, i_{\ast}i^{\ast}\Theta_{\PP^4}(1-d)) & \longrightarrow\\
         & \longrightarrow & H^2(\PP^4, \Theta_{\PP^4}(1-d) \otimes \mathcal{I}_X) & \longrightarrow & \cdots &  &  &
    \end{array}
\end{equation}
where $H^1(X, i^{\ast}\Theta_{\PP^4}(1-d)) \cong  H^1(\PP^4, i_{\ast}i^{\ast}\Theta_{\PP^4}(1-d))$.
For the first column of \eqref{TangSeq-2}, tensor the sequence (2) of Corollary \ref{K3-CI-IdShf} with the sheaf $\Theta_{\PP^4}(1-d)$ to obtain
\begin{equation*}
    0 \longrightarrow \Theta_{\PP^4}(-4-d) \longrightarrow \Theta_{\PP^4}(-1-d) \oplus \Theta_{\PP^4}(-2-d) \longrightarrow \Theta_{\PP^4}(1-d) \otimes \mathcal{I}_X \longrightarrow 0,
\end{equation*}
whose long exact cohomology sequence is
\begin{equation*}
    \begin{array}{c@{\hspace{3pt}}c@{\hspace{3pt}}c@{\hspace{3pt}}c@{\hspace{3pt}}c@{\hspace{3pt}}c@{\hspace{3pt}}c@{\hspace{3pt}}c}
         0 & \to & H^0(\PP^4, \Theta_{\PP^4}(-4-d)) & \to &  H^0(\PP^4, \Theta_{\PP^4}(-1-d)) \oplus H^0(\PP^4, \Theta_{\PP^4}(-2-d)) & \to & H^0(\PP^4, \Theta_{\PP^4}(1-d) \otimes \mathcal{I}_X) & \to\\
         & \to & H^1(\PP^4, \Theta_{\PP^4}(-4-d)) & \to &  H^1(\PP^4, \Theta_{\PP^4}(-1-d)) \oplus H^1(\PP^4, \Theta_{\PP^4}(-2-d)) & \to & H^1(\PP^4, \Theta_{\PP^4}(1-d) \otimes \mathcal{I}_X) & \to\\
         & \to & H^2(\PP^4, \Theta_{\PP^4}(-4-d)) & \to &  H^2(\PP^4, \Theta_{\PP^4}(-1-d)) \oplus H^2(\PP^4, \Theta_{\PP^4}(-2-d)) & \to & H^2(\PP^4, \Theta_{\PP^4}(1-d) \otimes \mathcal{I}_X) & \to\\
         & \to & H^3(\PP^4,\Theta_{\PP^4}(-4-d)) & \to &  \cdots &  &  & 
    \end{array}
\end{equation*}
Applying Serre duality and Lemma \ref{LemmBottFla}, we find:
\begin{itemize}
    \item 
    $H^1(\PP^4, \Theta_{\PP^4}(1-d))=0$,
    \item 
    $H^2(\PP^4, \Theta_{\PP^4}(-1-d))=H^2(\PP^4, \Theta_{\PP^4}(-2-d))=0$, and 
    \item $H^3(\PP^4, \Theta_{\PP^4}(-4-d))=0$ if and only if $d \neq 1$.
\end{itemize}
It follows that
\begin{equation}\label{dneq2}
    H^1(X,i^{\ast}\Theta_{\PP^4}(1-d)) = 0 \textup{ if and only if } d \neq 1,
\end{equation}
Combine \eqref{dim1ex2} and \eqref{dneq2} to conclude that
$H^1(X, \Theta_X \otimes i^{\ast}\mathcal{O}_{\PP^4}(1-d))=0$ for $d \geq 5$, which 
finishes the proof. 
\end{proof}
\end{theorem}

\begin{theorem}\label{Thm:Uniq-3}
    Let $ X \overset{ i }{\hookrightarrow} \PP^5 $ be a smooth complete intersection surface of type 
    $(2,2,2)$ and let $ d\geq 3 $ be an integer.
    Then, for $ d \geq 4 $, 
    any foliation $ \mathscr{F} = [s] \in \mathscr{F}ol( i^{\ast}\mathcal{O}_{\PP^5}(1-d), X ) $  
    with isolated singularities is uniquely determined by its singular scheme.     
    
\begin{proof}
Since $\mathcal{N}_{X/\PP^4} \cong i^{\ast}\mathcal{O}_{\PP^5}(2)^{\oplus 3}$, 
the sequence \eqref{NorCohSeqV2} becomes
\begin{equation*}
\begin{array}{c@{\hspace{3pt}}c@{\hspace{3pt}}c@{\hspace{3pt}}c@{\hspace{3pt}}c@{\hspace{3pt}}c@{\hspace{3pt}}c@{\hspace{3pt}}c}
     0 & \to & H^0(X, \Theta_X \otimes i^{\ast}\mathcal{O}_{\PP^5}(1-d)) & \to &  H^0(X, i^{\ast}\Theta_{\PP^5}(1-d)) & \to & H^0(X,i^{\ast}\mathcal{O}_{\PP^5}(3-d))^{\oplus 3} & \to \\
         & \to & H^1(X, \Theta_X \otimes i^{\ast}\mathcal{O}_{\PP^5}(1-d)) & \to &  H^1(X, i^{\ast}\Theta_{\PP^5}(1-d)) & \to & \cdots
    \end{array}
\end{equation*}
We compute the dimensions $h^0(X, i^{\ast}\mathcal{O}_{\PP^5}(3-d))$ and 
$h^1(X, i^{\ast}\Theta_{\PP^5}(1-d))$ separately:

To compute $h^0(X, i^{\ast}\mathcal{O}_{\PP^5}(3-d))$,
consider the sequence $(2)$ of Proposition \ref{PropHart}, with 
$ \mathcal{F} = \mathcal{O}_{\PP^5}(3-d) $. Apply Lemma \ref{aCM} to its 
associated cohomology sequence to obtain

\begin{equation}\label{Seq2-3}
 0  \longrightarrow  H^0(\PP^5, \mathcal{I}_X(3-d)) 
    \longrightarrow  H^0(\PP^5, \mathcal{O}_{\PP^5}(3-d))  \longrightarrow 
    H^0(\PP^5, i_{\ast}i^{\ast}\mathcal{O}_{\PP^5}(3-d))  \longrightarrow 0, 
\end{equation}
where $ H^0(\PP^5, i_{\ast}i^{\ast}\mathcal{O}_{\PP^5}(3-d)) 
\cong H^0(X, i^{\ast}\mathcal{O}_{\PP^5}(3-d)) $.

We claim that $ h^0(\PP^5, \mathcal{I}_X(3-d)) = 0 $: Indeed, since the homogeneous 
ideal of $ X $ is generated in degree $ 2 $, it follows that 
$ h^0(\PP^5, \mathcal{I}_X( t )) \geq 0 $ only for $ t \geq 2 $ (and $ 3 - d \leq 0 $).

Applying Lemma \ref{LemmBottFla} to \eqref{Seq2-3}, gives the result:

\begin{equation}\label{dim1ex3}
    3 \cdot h^0(X, i^{\ast}\mathcal{O}_{\PP^5}(3-d) )= 
    \begin{cases}
    3, & \text{if } d=3,\\
    0, & \text{if } d \geq 4.
    \end{cases}
\end{equation}

Next, to compute $h^1(X, i^{\ast}\Theta_{\PP^5}(1-d))$, consider the short exact sequence in $\PP^5$
$$0 \longrightarrow \Theta_{\PP^5}(1-d) \otimes \mathcal{I}_X \longrightarrow \Theta_{\PP^5}(1-d) \longrightarrow i_{\ast}i^{\ast}\Theta_{\PP^5}(1-d) \longrightarrow 0,$$
given by Proposition \ref{PropHart}. This sequence induces the long exact sequence
\begin{equation}\label{TangSeq-3}
    \begin{array}{c@{\hspace{5pt}}c@{\hspace{5pt}}c@{\hspace{5pt}}c@{\hspace{5pt}}c@{\hspace{5pt}}c@{\hspace{5pt}}c@{\hspace{5pt}}c}
         0 & \longrightarrow & H^0(\PP^5, \Theta_{\PP^5}(1-d) \otimes \mathcal{I}_X) & \longrightarrow &  H^0(\PP^5, \Theta_{\PP^5}(1-d)) & \longrightarrow 
         & H^0(\PP^5, i_{\ast}i^{\ast}\Theta_{\PP^5}(1-d)) & \longrightarrow\\
         & \longrightarrow & H^1(\PP^5, \Theta_{\PP^5}(1-d) \otimes \mathcal{I}_X) 
         & \longrightarrow &  H^1(\PP^5, \Theta_{\PP^5}(1-d)) & \longrightarrow 
         & H^1(\PP^5, i_{\ast}i^{\ast}\Theta_{\PP^5}(1-d)) & \longrightarrow\\
         & \longrightarrow & H^2(\PP^5, \Theta_{\PP^5}(1-d) \otimes \mathcal{I}_X) 
         & \longrightarrow & \cdots &  &  &
    \end{array}
\end{equation}
where $H^1(X, i^{\ast}\Theta_{\PP^5}(1-d)) \cong  H^1(\PP^5, i_{\ast}i^{\ast}\Theta_{\PP^5}(1-d))$. For the first column of \eqref{TangSeq-3}, tensor the sequence (3) of Corollary \ref{K3-CI-IdShf} 
with the sheaf $\Theta_{\PP^5}(1-d)$ to obtain
\begin{equation*}
    0 \longrightarrow \left(\rfrac{\mathcal{O}_{\PP^5}(-4)^{\oplus 3}}{\mathcal{O}_{\PP^5}(-6)}\right) \otimes \Theta_{\PP^5}(1-d) \longrightarrow \Theta_{\PP^5}(-1-d)^{\oplus 3} \longrightarrow \Theta_{\PP^5}(1-d) \otimes \mathcal{I}_X \longrightarrow 0,
\end{equation*}
whose long exact cohomology sequence is
\begin{equation}\label{Seq:LongEx32}
    \begin{array}{c@{\hspace{3pt}}c@{\hspace{3pt}}c@{\hspace{3pt}}c@{\hspace{3pt}}c@{\hspace{3pt}}c@{\hspace{3pt}}c@{\hspace{3pt}}c}
         0 & \to & H^0\left(\PP^5, \left(\rfrac{\mathcal{O}_{\PP^5}(-4)^{\oplus 3}}{\mathcal{O}_{\PP^5}(-6)}\right) \otimes \Theta_{\PP^5}(1-d)\right) & \to &  H^0(\PP^5, \Theta_{\PP^5}(-1-d))^{\oplus 3} & \to & H^0(\PP^5, \Theta_{\PP^5}(1-d) \otimes \mathcal{I}_X) & \to\\
         & \to & H^1\left(\PP^5, \left(\rfrac{\mathcal{O}_{\PP^5}(-4)^{\oplus 3}}{\mathcal{O}_{\PP^5}(-6)}\right) \otimes \Theta_{\PP^5}(1-d)\right) & \to &  H^1(\PP^5, \Theta_{\PP^5}(-1-d))^{\oplus 3} & \to & H^1(\PP^5, \Theta_{\PP^5}(1-d) \otimes \mathcal{I}_X) & \to\\
         & \to & H^2\left(\PP^5, \left(\rfrac{\mathcal{O}_{\PP^5}(-4)^{\oplus 3}}{\mathcal{O}_{\PP^5}(-6)}\right) \otimes \Theta_{\PP^5}(1-d)\right) & \to &  H^2(\PP^5, \Theta_{\PP^5}(-1-d))^{\oplus 3} & \to & H^2(\PP^5, \Theta_{\PP^5}(1-d) \otimes \mathcal{I}_X) & \to\\
         & \to & H^3\left(\PP^5, \left(\rfrac{\mathcal{O}_{\PP^5}(-4)^{\oplus 3}}{\mathcal{O}_{\PP^5}(-6)}\right) \otimes \Theta_{\PP^5}(1-d)\right) & \to &  \cdots &  &  & 
    \end{array}
\end{equation}
For the first column of \eqref{Seq:LongEx32}, tensor the sequence \eqref{New}
with the sheaf $\Theta_{\PP^5}$, whose long exact cohomology sequence is
\begin{equation*}
    \begin{array}{c@{\hspace{5pt}}c@{\hspace{5pt}}c@{\hspace{5pt}}c@{\hspace{5pt}}c@{\hspace{5pt}}c@{\hspace{5pt}}c@{\hspace{5pt}}c}
         0 & \to & H^0(\PP^5,\Theta_{\PP^5}(-5-d)) & \to &  H^0(\PP^5, \Theta_{\PP^5}(-3-d))^{\oplus 3} & \to & H^0\left(\PP^5,\left(\rfrac{\mathcal{O}_{\PP^5}(-4)^{\oplus 3}}{\mathcal{O}_{\PP^5}(-6)}\right) \otimes \Theta_{\PP^5}(1-d)\right) & \to\\
         & \to & H^1(\PP^5,\Theta_{\PP^5}(-5-d)) & \to &  H^1(\PP^5, \Theta_{\PP^5}(-3-d))^{\oplus 3} & \to & H^1\left(\PP^5,\left(\rfrac{\mathcal{O}_{\PP^5}(-4)^{\oplus 3}}{\mathcal{O}_{\PP^5}(-6)}\right) \otimes \Theta_{\PP^5}(1-d)\right) & \to\\
         & \to & H^2(\PP^5,\Theta_{\PP^5}(-5-d)) & \to &  H^2(\PP^5, \Theta_{\PP^5}(-3-d))^{\oplus 3} & \to & H^2\left(\PP^5,\left(\rfrac{\mathcal{O}_{\PP^5}(-4)^{\oplus 3}}{\mathcal{O}_{\PP^5}(-6)}\right) \otimes \Theta_{\PP^5}(1-d)\right) & \to\\
         & \to & H^3(\PP^5,\Theta_{\PP^5}(-5-d)) & \to &  H^3(\PP^5, \Theta_{\PP^5}(-3-d))^{\oplus 3} & \to & H^3\left(\PP^5,\left(\rfrac{\mathcal{O}_{\PP^5}(-4)^{\oplus 3}}{\mathcal{O}_{\PP^5}(-6)}\right) \otimes \Theta_{\PP^5}(1-d)\right) & \to\\
         & \to & H^4(\PP^5,\Theta_{\PP^5}(-5-d)) & \to &  \cdots & & &
    \end{array}
\end{equation*}
Applying Serre duality and Lemma \ref{LemmBottFla}, we find:
\begin{itemize}
    \item $H^1(\PP^5, \Theta_{\PP^5}(1-d))=0$,
    \item $H^2(\PP^5, \Theta_{\PP^5}(-1-d))=0$,
    \item $H^3(\PP^5, \Theta_{\PP^5}(-3-d))=0$,
    \item $H^4(\PP^5, \Theta_{\PP^5}(-5-d))=0$ if and only if $d \neq 1$.
\end{itemize}
It follows that
\begin{equation}\label{dneq3}
    H^1(X,i^{\ast}\Theta_{\PP^5}(1-d)) = 0 \textup{ if and only if } d \neq 1,
\end{equation}

Combine \eqref{dim1ex3} and \eqref{dneq3}, to conclude that
$H^1(X, \Theta_X \otimes i^{\ast}\mathcal{O}_{\PP^5}(1-d))=0$ for $d \geq 4$, which completes the proof. 
\end{proof}
\end{theorem}


\begin{thebibliography}{9}

\bibitem{Claudia}
C.R. Alcántara, Foliations on {CP}$^2$ of degree $d$ with a singular point with Milnor number $d^2+d+1$, \emph{Rev. Mat. Complut.} \textbf{31} (2018), 187–199.

\bibitem{ClaudiaRubi}
C.R. Alcántara and R. Pantaleón-Mondragón, Foliations on {CP}$^2$ with a unique singular point without invariant algebraic curves, \emph{Geom. Dedicata} \textbf{207} (2020), 193-200.

\bibitem{ClaudiaRonzon}
C.R. Alcántara and R. Ronzón-Lavie, Classification of Foliations on {CP}$^2$ of degree $3$ with Degenerate Singularities, \emph{J. Singul.} \textbf{14} (2016), 52-73.

\bibitem{A-C}
C. Araujo and M. Corrêa, On degeneracy schemes of maps of vector bundles and applications to
holomorphic foliations, \emph{Math. Z.} \textbf{276} (2014), 505-515.

\bibitem{Polarity}
A. Campillo and J. Olivares, Polarity with respect ot a foliation and {C}ayley-{B}acharach {T}heorems, \emph{J. Reine Angew. Math.} \textbf{534} (2001), 95-118.

\bibitem{CampilloOlivares}
A. Campillo and J. Olivares, On sections with isolated singularities of twisted bundles and applications to foliations by curves, \emph{Math. Res. Lett.} \textbf{10} (2003), 651-658.

\bibitem{C-F-N-V}
M. Corrêa Jr, A. Fernández-Pérez, G. Nonato Costa and R. Vidal Martins, Foliations by curves  with curves as singularities, \emph{Ann. Inst. Fourier Grenoble} \textbf{64} (2014), 1781-1805.

\bibitem{Cr-Es}
J. Cruz and E. Esteves, Regularity of subschemes invariant under {P}faff fields on projective spaces, \emph{Comment. Math. Helv.} \textbf{86} (2009), 947-965.

\bibitem{Friedman}
R. Friedman, \emph{Algebraic Surfaces and Holomorphic Vector Bundles}, Universitext, Springer New York, 2012.

\bibitem{RRAlgebra}
W. Fulton and S. Lang, \emph{Riemann-Roch Algebra}, Springer-Verlag, 1985.

\bibitem{Hirzebruch}
C. Galindo, F. Monserrat and J. Olivares, Foliations with isolated singularities on {H}irzebruch surfaces, \emph{Forum Math.} \textbf{33} (2021), 1471-1486.

\bibitem{G-P}
L. Giraldo and A.J. Pan-Collantes, On the singular scheme of codimension one holomorphic foliations in {P}$^3$, \emph{Int. J. Math.} \textbf{21} (2010), 843-858.

\bibitem{RuledSurfaces}
X. Gómez-Mont, Holomorphic foliations in ruled surfaces, \emph{Trans. Amer. Math. Soc.} \textbf{312} (1989), 179-201.

\bibitem{GMKempf}
X. Gómez-Mont and G. Kempf, Stability of meromorphic vector fields in projective spaces, \emph{Comm. Math. Helv.} \textbf{64} (1989), 462-473.

\bibitem{Hartshorne}
R. Hartshorne, \emph{Algebraic Geometry}, Graduate Texts in Mathematics, Springer, 1977.

\bibitem{Huybrechts}
D. Huybrechts, \emph{Lectures on K3 Surfaces}, Cambridge Studies in Advanced Mathematics, Cambridge University Press, 2016.

\bibitem{Okonek}
C. Okonek, M. Schneider and H. Spindler, \emph{Vector Bundles on Complex Projective Spaces: With an Appendix by S. I. Gelfand}, Modern Birkh{\"a}user Classics, Springer Basel, 2011.

\bibitem{MarcioInv}
M.G. Soares, The {P}oincaré problem for hypersurfaces invariant by one-dimensional foliations, \emph{Invent. Math.} \textbf{128} (1997), 495-500.

\bibitem{MarcioAnn}
M.G. Soares, Projective Varieties Invariant by One-Dimensional Foliations, \emph{Ann. of Math.} \textbf{152} (2000), 369-382.

\end{thebibliography}
\end{document}